\documentclass[a4paper]{article}
\usepackage{amsthm,amsmath,amssymb}
\newtheorem{teor}{Theorem}[section]
\newtheorem{defin}[teor]{Definition}
\newtheorem{lemm}[teor]{Lemma}
\newtheorem{osse}[teor]{Remark}
\newtheorem{prop}[teor]{Proposition}
\newtheorem{defi}[teor]{Definition}
\newtheorem{coro}[teor]{Corollary}
\newtheorem{prob}[teor]{Problem}

\newcommand{\bele}{\begin{lemm}\begin{sl}}
\newcommand{\enle}{\end{sl}\end{lemm}}
\newcommand{\bedef}{\begin{defi}\begin{sl}}
\newcommand{\eddef}{\end{sl}\end{defi}}
\newcommand{\bete}{\begin{teor}\begin{sl}}
\newcommand{\ente}{\end{sl}\end{teor}}
\newcommand{\beos}{\begin{osse}\begin{rm}}
\newcommand{\eddos}{\end{rm}\end{osse}}
\newcommand{\bepr}{\begin{prop}\begin{sl}}
\newcommand{\empr}{\end{sl}\end{prop}}
\newcommand{\bepro}{\begin{prob}\begin{rm}}
\newcommand{\empro}{\end{rm}\end{prob}}
\newcommand{\bede}{\begin{defin}\begin{sl}}
\newcommand{\edde}{\end{sl}\end{defin}}
\newcommand{\beco}{\begin{coro}\begin{sl}}
\newcommand{\enco}{\end{sl}\end{coro}}

\newcommand{\thspace}{\hspace{3mm}}

\newcommand{\quext}{\quad\text}
\newcommand{\qquext}{\qquad\text}
\newcommand{\de}{\partial}

\newcommand{\RR}{\mathbb{R}}
\newcommand{\NN}{\mathbb{N}}

\newcommand{\beeq}[1]{\begin{equation}\label{#1}}
\newcommand{\eddeq}{\end{equation}}

\newcommand{\beeqa}[1]{\begin{eqnarray}\label{#1}}
\newcommand{\eddeqa}{\end{eqnarray}}

\newcommand{\beal}[1]{\begin{align}\label{#1}}
\newcommand{\eddal}{\end{align}}

\newcommand{\bespl}[1]{\begin{split}\label{#1}}
\newcommand{\edspl}{\end{split}}

\newcommand{\bega}[1]{\begin{gather}\label{#1}}
\newcommand{\edga}{\end{gather}}

\newcommand{\beeqax}{\begin{eqnarray*}}
\newcommand{\eddeqax}{\end{eqnarray*}}

\def\qed{\ifmmode 
  \else \leavevmode\unskip\penalty9999 \hbox{}\nobreak\hfill
  \fi
  \quad\hbox{\hskip.5em\vrule width.4em height.6em depth.05em\hskip.1em}}
\def\endproofsym{\qed}
\newcommand{\dimbox}{\hbox{\hskip.5em\vrule width.4em height.6em depth.05em\hskip.1em}}
\renewenvironment{proof}[1][Proof]{\trivlist\item[\hskip\labelsep{\hskip0pt
    {\normalfont\scshape#1.}\hskip .321429\parindent}]\ignorespaces}
{\endproofsym\endtrivlist}
\def\endnobox{\def\endproofsym{}\end{proof}\def\endproofsym{\qed}}

\newcommand{\no}{\nonumber}

\newcommand{\beeqao}{\begin{eqnarray}\no}
\newcommand{\bealo}{\begin{align}\no}
\newcommand{\besplo}{\begin{split}\no}
\newcommand{\begao}{\begin{gather}\no}

\newcommand{\duav}[1]{\langle{#1}\rangle}

\newcommand{\cc}{{\mathfrak c}}

\newcommand{\perogni}{\forall\,}
\newcommand{\esiste}{\exists\,}

\newcommand{\io}{\int_\Omega}

\newcommand{\iTT}{\int_0^T}

\newcommand{\iTo}{\iTT\!\io}

\newcommand{\epsi}{\varepsilon}
\newcommand{\ee}{_{\varepsilon}}
\newcommand{\lla}{_{\lambda}}

\newcommand{\OO}{_{\Omega}}

\newcommand{\bn}{\boldsymbol{n}}
\newcommand{\bnn}{\boldsymbol{n}}

\newcommand{\dn}{\partial_{\bn}}

\newcommand{\lhs}{left hand side}
\newcommand{\rhs}{right hand side}

\DeclareMathOperator{\dive}{div}
\DeclareMathOperator{\deriv}{d}

\DeclareMathOperator{\Id}{Id}

\DeclareMathOperator{\entr}{entr}

\newcommand{\HUH}{H^1(0,T;H)}

\newcommand{\HUVp}{H^1(0,T;V')}

\newcommand{\LDH}{L^2(0,T;H)}
\newcommand{\LDV}{L^2(0,T;V)}
\newcommand{\LDVp}{L^2(0,T;V')}
\newcommand{\LIH}{L^\infty(0,T;H)}
\newcommand{\LIV}{L^\infty(0,T;V)}

\newcommand{\LDHD}{L^2(0,T;H^2(\Omega))}

\let\TeXchi\chi
\def\chi{{\setbox0 \hbox{\mathsurround0pt
$\TeXchi$}\hbox{\raise\dp0 \copy0 }}}

\newcommand{\kjk}{_k^{j_k}}

\newcommand{\zkk}{_{0,k}}
\newcommand{\jk}{^j_k}
\newcommand{\jkt}{^j_{k,t}}

\newcommand{\mmj}{_{m,j}}

\newcommand{\Wnj}{W_{n_j}}
\newcommand{\jjt}{_{j,t}}
\newcommand{\zzn}{_{0,n}}

\newcommand{\muciapo}{\widehat{\mu}}

\newcommand{\calH}{{\mathcal H}}

\newcommand{\calA}{{\mathcal A}}

\newcommand{\calE}{{\mathcal E}}
\newcommand{\calN}{{\mathcal N}}
\newcommand{\calEn}{{\mathcal E}_n}
\newcommand{\calM}{{\mathcal M}}

\newcommand{\calS}{{\mathcal S}}
\newcommand{\calSe}{{\mathcal S}_\epsi}

\newcommand{\calV}{{\mathcal V}}
\newcommand{\calY}{{\mathcal Y}}
\newcommand{\calW}{{\mathcal W}}
\newcommand{\calVm}{{\mathcal V}_m}
\newcommand{\calWm}{{\mathcal W}_m}
\newcommand{\dV}{\deriv_{\mathcal V}}

\newcommand{\dW}{\deriv_{\mathcal W}}

\newcommand{\calVe}{{\mathcal V}_\epsi}
\newcommand{\calMe}{{\mathcal M}_\epsi}

\newcommand{\Pee}{{\rm (P$_\epsi$)}}
\newcommand{\Pnee}{{\rm (P$_{n,\epsi}$)}}
\newcommand{\Pneen}{{\rm (P$_{n,\epsi_n}$)}}

\newcommand{\barO}{\overline{\Omega}}

\newcommand{\dit}{\deriv\!t}
\newcommand{\dis}{\deriv\!s}
\newcommand{\dix}{\deriv\!x}
\newcommand{\dir}{\deriv\!r}

\newcommand{\ddt}{\frac{\deriv\!{}}{\dit}}

\newcommand{\calAen}{{\calA}_{\entr}}
\newcommand{\calSen}{{\calS}_{\entr}}

\numberwithin{equation}{section}
\begin{document}

\title{Global attractors 
  for Cahn-Hilliard equations with non constant mobility}
\author{Giulio Schimperna\\
{\sl Dipartimento di Matematica, Universit\`a di Pavia}\\
{\sl Via Ferrata, 1}\\
{\sl I-27100 Pavia, Italy}\\
{\tt giusch04@unipv.it}\\
}
\maketitle
\begin{abstract}
 We address, in a three-dimensional spatial setting, 
 both the viscous and the standard Cahn-Hilliard equation
 with a nonconstant mobility coefficient. As it was shown in
 J.W.~Barrett and J.W.~Blowey, Math.\ Comp., 68 (1999),
 487--517, one cannot expect uniqueness of the solution
 to the related initial and boundary value problems. 
 Nevertheless, referring to J.~Ball's theory of generalized
 semiflows, we are able to prove existence of 
 compact quasi-invariant global attractors for the associated
 dynamical processes settled in the natural ``finite
 energy'' space. A key point in the proof is a 
 careful use of the energy equality, combined with 
 the derivation of a ``local compactness'' estimate for systems 
 with supercritical nonlinearities, which may have 
 an independent interest. Under growth restrictions
 on the configuration potential, we also show existence
 of a compact global attractor for the semiflow generated
 by the (weaker) solutions to the nonviscous equation
 characterized by a ``finite entropy'' condition.
\end{abstract}


\vspace{2mm}

{\bf AMS (MOS) subject clas\-si\-fi\-ca\-tion:}\thspace
35K55, 35B41, 35B45.



%
%


\section{Introduction}
\label{secintro}

In this note we address the initial and (homogeneous Neumann)
boundary value problem for the equation
\begin{equation}\label{CHin}
  u_t-\dive\big(b(u)\nabla (\epsi u_t-\Delta u+W'(u)+f)\big),
\end{equation}
which is settled in a smooth and bounded domain $\Omega\subset\RR^3$
and corresponds for $\epsi=0$ to the standard,
and for $\epsi>0$ to the {\sl viscous}, Cahn-Hilliard
equation with {\sl nonconstant}\/ mobility function
$b(\cdot)$. In particular, $b$ is assumed to depend 
on $u$ in a globally Lipschitz way and 
is not allowed to degenerate. In the relation above,
$W$ is a possibly nonconvex configuration
potential and $f$ a source 
which is included in view of possible applications
to conserved phase field models
(where $u$ is an order parameter
and $f$ represents a coupling term
depending on the temperature, 
see, e.g.~\cite{Ca}). Relation \eqref{CHin}
is complemented by homogeneous Neumann boundary
conditions both for $u$ and the {\sl chemical potential}\/
$w:=(\epsi u_t-\Delta u+W'(u)+f)$.

The dependence of the mobility on the variable $u$ 
is very relevant for physical applications. Actually,
as $u$ represents the density of one component in a binary
alloy, one expects that the diffusion of mass is 
influenced by the actual configuration, i.e., by the 
value of $u$. In fact, it is just due to difficulties
arising in the analysis of \eqref{CHin}
that, in the mathematical literature, $b$ has been 
generally replaced by a constant function. 

The most relevant work devoted to the mathematical 
study of \eqref{CHin} for nonconstant (but nondegenerate)
$b$ is \cite{BB},
where, for zero source $f$ and no viscosity 
(i.e., $\epsi=0$),
existence and uniqueness of the solution, 
together with additional regularity properties, are proved in space 
dimensions 1 and 2. On the contrary,
in the three dimensional case, only existence of 
a weak solution is shown (uniqueness would
hold in 3D for a class of more regular solutions, but 
the authors cannot prove this further regularity).
The results of \cite{BB}, which are also complemented
by numerical investigations, are very sharp and it seems
rather difficult to fill the regularity gap which 
prevents from having well posedness in 3D. Actually,
some more recent work \cite{LQY} has been devoted
to improve the regularity of solutions, but
still only in the 2D case.

Here, we aim to analyze, referring just to 
the 3D setting, the long time behavior
of \eqref{CHin} from the point of view of global 
attractors and considering both the viscous
and the nonviscous case. Due to the quoted difficulties,
this analysis is far from being trivial.
Actually, the use of more or less standard 
tools seems possible only for the viscous equation and if the 
potential $W$ has a controlled growth at $\infty$
(cf.~\eqref{growW} below). Indeed, in this case,
uniqueness holds at least for $t>0$ (for $t\ge0$
if the initial datum is more regular) and we have   
uniform regularization properties. Instead,
if $\epsi=0$ and/or we are in the 
(physically relevant) situation of 
fastly growing or even {\sl singular}\/ 
(i.e., uniformly taking the value $+\infty$ outside
a bouded interval, cf.~\eqref{singW} below) potentials,
we then have to proceed much more carefully. 

Actually, in such a framework, existence of solutions 
and dissipativity of the process 
are still  easy to show, but then we have to face 
the following three main difficulties:
%
%
\begin{enumerate}
\item
 We have no uniqueness result. Thus, we have to
 refer to some machinery which is suitable for dealing
 with problems with lack of uniqueness. Among the various 
 possible choices (we quote in particular the alternative
 possibility to work the in space of trajectories,
 cf., e.g., \cite{ChV,Se}), we decided to refer 
 to J.~Ball's theory of {\sl generalized semiflows}\/
 \cite{Ba1,Ba2} which has the advantage 
 of being very close to the standard physical 
 interpretation. Namely, the system still gives rise to a 
 dynamical process settled in a 
 phase space $\calV$ of states (rather 
 than, for instance, of trajectories).
 To be more precise, due to point (ii)~below, 
 a further generalization of
 Ball's approach, recently devised in \cite{RSS}
 (see also \cite{MV}), will be used.
\item
 Not all the estimates we perform can be rigorously
 carried out in the regularity framework which appears
 to be the natural one for~\eqref{CHin}.
 Namely, one has to proceed through approximation and 
 passage to the limit. However, due to lack of uniqueness,
 it is not obvious whether all the solutions with  
 the natural regularity can be reached by
 the approximation procedure. Thus, we have to restrict
 ourselves to consider solutions which are limit
 of more regular sequences for which the estimates 
 can be rigorously shown. This has a consequence on 
 the structure of the global attractor, which turns out
 to be only {\sl quasi-invariant}\/ rather than
 {\sl fully invariant}\/ as in the standard cases
 (cf.~\cite[Def.~2.8]{RSS} and Remark~\ref{noH3} below 
 for more details on this point).
\item
 Finally, despite the strictly parabolic character
 of the system, we cannot prove any uniform 
 in time regularization property
 of solutions (which, by the way, would also lead
 to uniqueness). For this reason, we have 
 to get the asymptotic compactness of the process 
 through a different and rather nonstandard procedure,
 which in our opinion can have an independent interest
 and might be applied to other systems with supercritical
 or fastly growing nonlinearities. Actually, we combine 
 the use of the energy equality, which is a consequence
 of the variational structure of \eqref{CHin}
 and is satisfied by all solutions in our regularity class, 
 with a ``locally uniform'' regularization property.
 Namely, we can show that there exist a set $K_0$,
 compact in the phase space $\calV$, 
 and a number $\delta>0$, both independent
 of the initial data, such that all admissible solutions 
 $u=u(t)$ starting from a given set $B$ bounded in $\calV$,
 after some $T_0>0$ depending only on the radius
 of $B$ in $\calV$ satisfy that
 \beeq{keypropintro}
   \perogni t\ge T_0,~~\esiste \tau=\tau(t)\in[0,3/2]:~~
    \perogni s\in [t+\tau,t+\tau+\delta],~~
    u(s)\in K_0.
 \end{equation}
 Such a property (combined with the energy equality)
 turns out to imply the asymptotic compactness of 
 the semiflow. Unfortunately, we are not able to show
 \eqref{keypropintro} for all ``singular potentials'' $W$
 (i.e.~those being $+\infty$ outside a bounded interval,
 here normalized to $(-1,1)$ for convenience), 
 but only for those of a subclass 
 (introduced in \cite[(H5)]{SSS1} and 
 called here of ``separating'' potentials, see 
 Def.~\ref{defsepW} below), which turn out to explode 
 sufficiently fast in proximity of $\pm1$.
 In particular, this class does {\sl not}\/ 
 seem to include the
 {\sl logarithmic}\/ potential
 \beeq{logpot}
   W(r)=(1+r)\log(1+r)
    +(1-r)\log(1-r)-\frac\lambda2r^2,
    \qquad r\in(-1,1),
 \end{equation}
 where $\lambda$ is a positive parameter,  
 relevant in concrete physical situations.
 We also remark that, due to the nonuniform character
 of \eqref{keypropintro}, the resulting global
 attractor will be compact in $\calV$, but not
 necessarily bounded in a ``better'' space.
\end{enumerate}
The procedure sketched above can be applied both
for $\epsi>0$ and for $\epsi=0$. Of course, if $\epsi>0$ and
$W$ has a controlled growth at infinity, we have
uniform regularization and, at least, unique continuation 
of solutions. Thus, the global attractor
can be intended in the framework of the 
standard theory for single-valued semigroups
(cf., e.g., \cite{Te}). Moreover, in this case one
could prove with only technical difficulties
further regularity properties of the attractor.

Finally, in the nonviscous case $\epsi=0$
we can also prove existence of a  
compact set in $\calV$ which uniformly attracts
{\sl less regular}\/ solutions, namely
those taking values in a larger
phase space $\calH$, 
which we call of ``finite entropy''. Actually,
by approximation, existence in this class 
(which was not considered in \cite{BB}) 
can be proved by means of an estimate
of entropy type (cf., e.g., \cite{DPGG}), 
for whose validity, however, a growth 
restriction on $W$ seems essential (excluding
from this result any singular potential).
Moreover, the entropy estimate turns out to
have a dissipative character, yielding
existence of an absorbing set in $\calH$.
Then, we prove that from any initial 
datum $u_0\in \calH$ there 
starts {\sl at least}\/ one solution $u$,
which, for $t>0$, lies in the energy
space $\calV$, which is compactly
embedded in $\calH$.
Clearly, this implies existence of a global attractor
bounded in the ``better'' space $\calV$. Of
course, also in this case, we have no uniqueness
and are still forced to use the ``generalized 
semiflows'' machinery described above.

The rest of this paper is organized as
follows. In the next Section~\ref{secmain}, after
recalling some preliminary material, 
we present our hypotheses and state our
main results, with the exception of those
related to entropy solutions. The related proofs
are given in the subsequent Section~\ref{secvisc}.
Finally, Section~\ref{secentropy} is devoted to 
the analysis of solutions in the finite entropy
class.


\section{Notations and main results}
\label{secmain}

Let $\Omega\subset\RR^3$ be a smooth bounded domain.
Let us set $H:=L^{2}(\Omega)$ and denote by
$(\cdot,\cdot)$ the scalar product in $H$ and by
$\|\cdot\|$ the related norm. The same symbols are
used also to note $H^3$ and its scalar product
and norm. The symbol $\|\cdot\|_{X}$ 
will indicate the norm in the generic Banach space $X$. 
Let us also assume that
\begin{equation} \label{hpb} 
  b\in W^{1,\infty}(\RR;\RR),\qquad
   \esiste \alpha,\mu>0:~~
   \alpha\le b(r)\le\mu~~\perogni r\in\RR.
\end{equation}
Then, we set $V:=H^1(\Omega)$, endowed with its standard
scalar product and norm. Letting $u:\Omega\to\RR$ be 
a {\sl measurable}\/ function, we 
introduce the couple of elliptic
operators $B,B_u:V\to V'$ (where $V'$ is the topological
dual of $V$), respectively given by
\beeq{defiBu}
  \langle Bv,z\rangle=\io\nabla v\cdot\nabla z,
   \qquad
  \langle B_u v,z\rangle=\io b(u)\nabla v\cdot\nabla z,
\end{equation}
the notation $\langle\cdot,\cdot\rangle$
standing for the duality between $V'$ and $V$. 
Then, we clearly have
\begin{equation} \label{coerc} 
  \langle B_u v,z\rangle\le \mu\|v\|_V\|z\|_V,\qquad
   \langle B_u v,v\rangle\ge \alpha\|\nabla v\|^2
\end{equation}
for all $u,v,z$ as before. If $u$ additionally
depends on time (i.e.~it is a measurable function
defined on $\Omega\times(0,T)$ for some $T>0$),
then $B_u$ is naturally extended to time dependent
functions. Namely, $B_u:L^2(0,T;V)\to L^2(0,T;V')$
is still a continuous and coercive operator
given by
\begin{equation}\label{But}
  \iTT\langle B_u v,z\rangle:=
   \iTo  b(u(x,t))\nabla v(x,t)\cdot\nabla z(x,t)\,\dix\,\dit.
\end{equation}
We shall adopt in the sequel the convention of
writing
\beeq{defivoo}
  \zeta\OO:=|\Omega|^{-1}\langle \zeta,1 \rangle
\end{equation}
for $\zeta\in V'$, where $|\Omega|$ stands for
the Lebesgue measure of $\Omega$. Let us also set
\beeq{defiV0}
  V_0':=\{\zeta\in V':\zeta\OO=0\},\quad
   H_0:=H\cap V_0',\quad
   V_0:=V\cap V_0
\end{equation}
and observe that, in general, if $v\in V$ and 
$\zeta\in V'$, then
\beeq{medie1}
  \langle \zeta-\zeta\OO,v \rangle
   = \langle \zeta-\zeta\OO,v-v\OO \rangle
   = \langle \zeta,v-v\OO \rangle.
\end{equation}
Moreover, if $u$ is as in \eqref{coerc},
then $B_u$ is bijective from 
$V_0$ to $V_0'$, so that we can define its 
inverse $\calN_u$, which fulfills,
for all $v\in V$, $\zeta\in V'$,
\beeq{medie2}
  \langle B_u v,\calN_u(\zeta-\zeta\OO) \rangle
   = \langle \zeta-\zeta\OO,v \rangle.
\end{equation}
Let us now come to the assumptions on the 
potential $W$. We let $I$ be an 
open interval of $\RR$ containing $0$
(possibly unbounded or even coinciding
with the whole real line),
$\lambda>0$, $c_W\ge 0$, and assume that
\beal{W}
  & W\in C^{2}(I;\RR),\quad W'(0)=0,\\
 \label{Wnew}
  & W(r)\ge 3\lambda r^2-c_W\quad\perogni r\in I,\\
 \label{coercW}
  & W''(r)\ge -\lambda\quad
   \perogni r\in I.
\end{align}
Assumption \eqref{Wnew} states
that the growth rate of $W$ for large values of $r$ 
is sufficiently fast to compensate its 
possible nonconvexity \eqref{coercW}
near 0. Of course, \eqref{Wnew} holds automatically 
whenever, for large $r$, it is $W(r)\sim \eta |r|^q$
for some $\eta>0$ and $q>2$.
Additionally, we shall assume either of properties
\eqref{growW}, \eqref{singW} below. The first 
is a {\sl controlled growth}\/ condition
(the choice of $p\in[2,\infty]$ will be 
made precise later):
\beeq{growW}
  I=\RR,\qquad \esiste K_W>0:~~
   W''(r)\le K_W (1+|r|^{p-2})~~
   \perogni r\in\RR,
\end{equation}
Of course, the larger is $p$, the weaker
is \eqref{growW} and, conventionally, we assume
that for $p=\infty$, \eqref{growW} just means $I=\RR$. 
The second condition identifies the so-called
{\sl singular}\/ potentials:
\beeq{singW}
  I=(-1,1),\qquad \lim_{|r|\to 1^-}
   W'(r)r=+\infty.
\end{equation}
In particular, in case
\eqref{singW} holds, then \eqref{Wnew}
is an immediate consequence of its.
Note that the domain $I$ of 
$W$ has been normalized to $(-1,1)$ just
for the sake of simplicity. 
We also let
\begin{equation}\label{regof}
  f\in H.
\end{equation} 
Let us now introduce
the {\sl energy}\/ of the system 
(possibly taking the value $+\infty$ for some 
$v$) as
\beeq{defiE}
  \calE(v):=\io\Big(\frac{|\nabla v|^2}2+W(v)+fv\Big),
   \quext{for }\,v\in V.
\end{equation}
Then, we define the space of data of {\sl finite energy}\/
as
\beeq{deficalV}
  \calV:=\big\{v\in V:~W(v)\in L^1(\Omega)\big\}.
\end{equation}
By continuity of the embedding $V\subset L^6(\Omega)$
it is clear that, if \eqref{growW} holds with $p\le 6$,
then it is in fact $\calV=V$. Otherwise, $\calV$ can be a 
proper subset of $V$.

We also set $\beta(r):=W'(r)+\lambda r$ for $r\in I$. On
account of \eqref{coercW}, $\beta$ is a monotone function
which will be sometimes identified with a {\sl maximal monotone}\/
operator from $H$ to itself (note that the maximality of
$\beta$ follows from the second condition in \eqref{singW}
if $W$ is singular). Then, we define 
\beeq{deficalW}
  \calW:=\big\{v\in H^2(\Omega):~\dn v=0~\text{on }\de\Omega,
   ~\beta(v)\in L^2(\Omega)\big\}.
\end{equation}
The set $\calW$ is nothing else 
than the domain (in $H$) of the {\sl subdifferential}\/ 
$\de \calE(v)$, where $\calE$ is
now seen as a (bounded from below) functional on $H$.
We can also introduce a metric structure on $\calV$
by setting
\beeq{defidV}
  \dV(v,z):=\|v-z\|
   +\|W(v)-W(z)\|_{L^1(\Omega)}
  \quad\perogni v,z\in \calV.
\end{equation}
Proceeding as in \cite[Lemma~3.8]{RS}, one can easily show
that $\calV$ is a complete metric space with the distance
$\dV$. Of course, the contribution of the second
term in the \rhs\ above is redundant, and could be
omitted so that $\calV=V$, in case \eqref{growW} holds
with $p\le 6$.
Analogously, $\calW$ is endowed with the distance
\beeq{defidW}
  \dW(v,z):=\|v-z\|_{H^2(\Omega)}
   +\|\beta(v)-\beta(z)\|
  \quad\perogni v,z\in \calW,
\end{equation}
where the second term on the \rhs\ is included
only in case \eqref{singW} holds; otherwise, it can be
omitted. It is clear that also $\calW$ is
a complete metric space.
%
%
Assuming \eqref{coercW} (respectively, \eqref{singW}),
we shall take $m>0$ (respectively, $m\in(0,1)$) and consider
the metric-closed subset of $\calV$ given by
\beeq{deficalVm}
  \calVm:=\big\{v\in \calV:~|v\OO|\le m\big\}.
\end{equation}
We also define, analogously, a closed subset $\calWm$
of $\calW$. The proof of the following result, which collects
further properties of $\calV$, $\calW$, and of the
energy $\calE$, can be performed by standard semicontinuity 
arguments (cf.~also \cite[Lemma~4.2]{RS}),
and it is thus omitted.
\bele\label{lemmaWcompV}
 The functional $\calE$ is sequentially weakly lower 
 semicontinuous on $V$. Moreover, $\calW$ is compactly embedded
 into $\calV$, namely any bounded sequence in $\calW$ admits
 a subsequence converging in $\calV$. 
 Finally, the convergence $v_n\to v$ in $\calV$ 
 is equivalent to the coupling of 
 \beeq{equiconv}
   v_n \to v\quext{weakly in }\,V 
    \qquext{and}\quad \limsup_{n\nearrow\infty}\calE(v_n)\le\calE(v).\dimbox
 \end{equation}
\enle
For $m$ as above, the initial datum $u_0$ is then chosen such that
\beeq{regou0}
  u_0\in\calVm.
\end{equation}
We are now ready to introduce our first notions of solutions 
to~\eqref{CHin}. We shall treat the viscous
($\epsi>0$) and the ``standard'' ($\epsi=0$) equation
altogether. 
\bede\label{defisolepsi}
 We call a (global) {\rm energy solution} to\/ {\rm Problem~\Pee}
 if $\epsi>0$ (respectively, to\/ {\rm Problem~(P$_0$)} if $\epsi=0$)
 one function $u:\Omega\times(0,\infty)\to \RR$ 
 such that, for all $T>0$, the regularity properties
 \beal{regou} 
    & u\in \HUVp\cap \LIV\cap \LDHD, \qquad \epsi^{1/2}u\in\HUH,\\
  \label{regow} 
    & W'(u)\in \LDH, \qquad  w\in \LDV
 \end{align}
 are fulfilled, and $u$ satisfies, in the space $V'$
 and for almost all times
 in $(0,\infty)$, the equations
 \beal{eqne}
   & u_t+B_u w=0,\\
  \label{eqne2}
   & w=\epsi u_t+Bu+W'(u)+f
 \end{align}
 and, a.e.~in $\Omega$, the initial condition
 \beeq{iniz}
   u|_{t=0}=u_0.
 \end{equation}
 Moreover, we say that an energy solution is\/ 
 {\rm regularizing} if the properties
 \beeq{regou2}
   u\in L^\infty(\tau,T;H^2(\Omega)),\quad
   u_t\in L^2(\tau,T;V),\quad
    \beta(u)\in L^\infty(\tau,T;H),
 \end{equation}
 hold for all $\tau>0$, $T\ge\tau$.
\edde
In the above statement, \eqref{CHin} has been split,
for convenience, as a system of the two
equations \eqref{eqne}--\eqref{eqne2}.
Testing \eqref{eqne} by 1, one
immediately sees that, for any energy solution 
to~\Pee\ or to (P$_0$)
corresponding to the initial datum $u_0$, it is
\beeq{consu}
  (u(t))\OO=(u_0)\OO=:u\OO
   \quad\perogni t\ge 0.
\end{equation}
Thus, $\calVm$ can be used as a {\sl phase space}\/
for the dynamical processes associated to Problems~\Pee, (P$_0$).

Let us now come to mathematical results, and we start by establishing
existence and, conditionally, uniqueness. 
\bete\label{teoesivisco}
 Let\/ \eqref{hpb}, \eqref{W}--\eqref{coercW},
 either\/ \eqref{growW} with $p=\infty$ 
 or\/ \eqref{singW}, \eqref{regof}
 and\/ \eqref{regou0} hold, and let $\epsi\ge 0$.
 Then, there  exists at least one energy solution $u$
 to\/ {\rm Problem~\Pee} (if $\epsi=0$,
 to {\rm Problem~(P$_0$)}\/).
 Moreover, if $u_1,u_2$ are a pair of energy solutions
 either to~\Pee\ or to~{\rm (P$_0$)}, satisfying,
 for some $\tau\ge 0$, property\/ \eqref{regou2}
 and such that $u_1(\tau)=u_2(\tau)$, then $u_1\equiv u_2$
 on $[\tau,\infty)$.
 Finally, only in case\/ $\epsi>0$, and if 
 \eqref{growW} holds with $p\le6$,
 then\/ {\rm Problem~\Pee} admits at least one\/
 {\rm regularizing} solution $u$ which, if $u_0\in\calW$,
 satisfies~\eqref{regou2} also for $\tau=0$.
\ente
The proofs of the above Theorem, and of the ones
which follow, are all posponed to the next Section.
Note that the existence part of the statement
above, at least for $\epsi=0$,
follows more or less
the lines of \cite[Thm.~2.2]{BB}, so that
we do not claim originality here.
Note also that, if $\epsi>0$ and 
\eqref{growW} holds with 
$p\le 6$, uniqueness is satisfied
starting from $\tau=0$ if $u_0\in\calW$,
and from any $\tau>0$ if $u_0\in\calV$
(namely, we have {\sl unique continuation}\/ of
trajectories). Instead, the uniqueness part might 
be vacuous (because we cannot prove existence 
of regularizing solutions) in all other cases
(in particular, if it is $\epsi=0$).

Theorem~\ref{teoesivisco} will be proved
by working on an approximate statement 
that we now introduce.
First of all, we replace $W$ by a
regularized potential $W_n$, 
with $n$ intended to go to $\infty$
in the limit, constructed this way. 
Recalling that $\beta=(W'+\lambda\Id)$ is 
monotone by \eqref{coercW}, we note 
as $\beta_n$ its Yosida approximation of index $n^{-1}$. 
Next, we define $W_n':=\beta_n-\lambda\Id$.
Then, $W_n'$ is (globally in $\RR$) Lipschitz continuous 
(the Lipschitz constant of course depending on $n$)
and it tends to $W'$ in the sense of {\sl G-convergence}\/ 
(see, e.g., \cite[Chap.~3]{At}). Moreover, 
defining $W_n$ by integration and choosing appropriately
the integration constant, one has that
\eqref{W}, \eqref{coercW} (and possibly 
\eqref{growW}) still hold for $W_n$, uniformly in $n$. 
Moreover, $W_n(r)\le W(r)$ for all $n\in\NN$,
$r\in I$, and, in place of \eqref{Wnew}, 
there holds, at least for $n$ sufficiently large,
\beeq{Wnewappr}
  W_n(r)\ge 2\lambda r^2-c_W\quad\perogni n\in\NN,~r\in \RR.
\end{equation}
Next, we replace $u_0\in\calV_m$
by a regularizing sequence 
$\{u\zzn\}\subset H^2(\Omega)\cap \calV_m$ tending
to $u_0$ in $V$ for $n\nearrow\infty$.
Namely, we define $u\zzn$ as the unique solution
to the elliptic problem
\beeq{uzn}
  u\zzn\in H^2(\Omega),\qquad
   \frac1n Bu\zzn+u\zzn=u_0.
\end{equation}
Then, it is well known \cite{Li2} that $(u\zzn)\subset H^3(\Omega)$,
$u\zzn\to u_0$ strongly in $V$, and
\beeq{uzn2}
  \|u\zzn\|_V\le \|u_0\|_V, \qquad
   \|u\zzn-u_0\|\le n^{-1/2}\|u_0\|_V
   \quad\perogni n\in\NN,
\end{equation}
so that, by standard properties of subdifferentials
and using \eqref{uzn}, we also have
\bealo
  \io W_n(u\zzn)
   & \le \io W_n(u_0)+\big(\beta_n(u\zzn),u\zzn-u_0\big)
    - \frac\lambda2\|u\zzn\|^2 + \frac\lambda2\|u_0\|^2\\
 \label{Wuzn}
  & \le \io W_n(u_0)
    - \frac\lambda2\|u\zzn\|^2 + \frac\lambda2\|u_0\|^2
    \le \io W(u_0) + \sigma_n,
\end{align}
where $\sigma_n$ goes to 0 as $n\nearrow\infty$. 
Then, if $\epsi>0$, the replacements of $W$ with $W_n$ and
of $u_0$ with $u\zzn$ give rise to a new Problem~\Pnee. 
If $\epsi=0$, we additionally take $\epsi=\epsi_n>0$ 
in \eqref{eqne2}, where $(\epsi_n)\subset(0,1)$ is some sequence
going to $0$ as $n\nearrow\infty$
(e.g., $\epsi_n=n^{-1}$), and we get a 
Problem we call~\Pneen. It is clear that, at least formally,
Problems \Pnee\ and \Pneen\ tend, as $n\nearrow\infty$,
to~\Pee\ and (P$_0$), respectively.
By a standard application of the Faedo-Galerkin method 
(cf.~\cite{BB} for the details), one can easily
show that, for every $n\in\NN$, each of Problems~\Pnee, \Pneen\
admits one and only one (global in time) solution 
(in both cases we note it by $u_n$).
Actually, both the uniqueness property and the global
character of the solutions are not directly
guaranteed by the Faedo-Galerkin method, 
but they can be shown proceeding along the lines
of the next Section (we thus omit the details). 
Moreover, setting for $T>0$
\beeq{YT}
   \calY_T:=H^2(0,T;H)\cap W^{1,\infty}(0,T;V)\cap
       H^1(0,T;H^2(\Omega)),
\end{equation}
it can be proved with only technical difficulties that,
for all $n\in\NN$,
\beeq{regoun}
  u_n\in \calY_T \qquad \perogni T>0.
\end{equation}
Let us then come to the long-time issue. To begin with, 
we need some preliminary work, starting with a simple
property satisfied by all solutions.
\bele\label{enerbuone}
 Let $u$ be an energy solution either to\/ {\rm Problem~\Pee} 
 or to\/ {\rm Problem~(P$_0$)} and let $T>0$.
 Then, $u(t)\in \calV_m$ for\/ {\rm all} (not just a.e.)
 $t\in[0,T]$. Moreover, for\/
 {\rm all} $t,t_1,t_2\ge 0$, $u$ satisfies the\/
 {\rm energy equality}
 \beeq{energyeq}
  \calE(u(t_2))- \calE(u(t_1))
   +\int_{t_1}^{t_2}\io \big(b(u)|\nabla w|^2\big)
   +\epsi\int_{t_1}^{t_2}\|u_t\|^2=0
 \end{equation} 
 and the\/ {\rm dissipativity estimate}
 \beeq{dissepsi}
  \calE(u(t))\le\calE(u_0)e^{-\kappa t}+C_0,
 \end{equation}
 where $\kappa,C_0>0$ are computable constants,
 independent both of the initial
 data and of $\epsi$ (they can depend on $m$,
 cf.~\eqref{deficalVm}, instead).
\enle
However, to define the dynamical processes associated 
to~\Pee\ and (P$_0$), it seems necessary to restrict the classes
of admissible solutions. Actually, while the above Lemma
shows uniform dissipativity for {\sl all}\/ energy solutions,
as we look for some form of parabolic regularization in
time, we readily realize that energy solutions 
need not be smooth enough to prove rigorously 
sharper estimates. For this reason,
we ``essentially'' (cf.~Remark~\ref{nonproprioPn} below)
have to work on Problem~\Pnee\ (or on
Problem~\Pneen) and then take the limit $n\nearrow\infty$.
Let us then introduce a useful notion of convergence: given
$T>0$, we say that a sequence $(u_j)$ of functions 
from $\Omega\times(0,T)$ to $\RR$ tends to $u$ 
{\sl weakly in $\calV_T$} if the following 
two properties hold:
\beal{conve1}
   & u_j\to u \quext{weakly in }\,\HUVp\cap\LDHD 
    \quext{(also in $\,\HUH\,$, if $\,\epsi>0$)},\\
  \label{conve2}
   & u_j\to u \quext{weakly star in }\,\LIV.
\end{align}
We can thus introduce the class of solutions for which
we shall prove existence of the global attractor.
\bede\label{defisollimi}
 Let $u$ be an {\rm energy solution} to\/ {\rm Problem~\Pee}
 (respectively, to\/ {\rm Problem~(P$_0$)}\/). We say that 
 $u$ is\/ {\rm limiting} if there exists a sequence
 $(u_j)$, with $(u_j)\subset \calY_T$ for all
 $T>0$, such that, as $j\nearrow\infty$, $u_j$ tends to 
 $u$ weakly in $\calV_T$ for all $T>0$; moreover, for 
 each $j\in\NN$, there exists $n_j\in \NN$,
 with $n_j\nearrow\infty$ for $j\nearrow\infty$
 (if $\epsi=0$ we also ask existence of 
 $(0,1)\ni\epsi_j\to0$),
 such that $u_j$ solves, in the usual sense,
 \beal{eqnej}
   & u_{j,t}+B_{u_j} w_j=0,\\
  \label{eqne2j}
   & w_j=\epsi_j u_{j,t}+Bu_j+W_{n_j}'(u_j)+f
 \end{align}
 (with $\epsi_j\equiv\epsi$ if $\epsi>0$)
 and, finally, there hold
 \bega{conve3}
    W_{n_j}'(u_j)\to W'(u) \quext{weakly in }\,\LDH,
      \qquad\perogni T>0,\\
  \label{conve4} 
   \esiste\sigma_j\searrow 0:\quad
    \calE_{n_j}(u_j(0)) \le \calE(u_0) + \sigma_j.
 \end{gather}
 Here, for $n\in\NN$, the approximate energy $\calE_n$ is defined
 as in \eqref{defiE}, but with $W_n$ replacing $W$.
\edde
\beos\label{nonproprioPn}
 In the proof of Theorem~\ref{teoesivisco} we shall show
 that the class of limiting solutions is not empty
 by passing to the limit  $n\nearrow\infty$ in \Pnee\ 
 (or \Pneen). In particular, we then have 
 that \eqref{conve4} is satisfied thanks to
 \eqref{uzn2} and \eqref{Wuzn}. However, this 
 natural procedure seems not ``robust'' enough 
 (especially from the viewpoint of taking subsequences)
 to guarantee that the resulting set of solutions
 satisfies the desired semiflow properties 
 (in particular, (H4) in Definition~\ref{defilimsem} below).
 This is the reason why in Definition~\ref{defisollimi} 
 above we are forced to generalize a bit the method
 (which, as a further consequence, possibly enlarges the class
 of limiting solutions).
\eddos
%
%
\beos\label{escluse}
 Clearly, due to nonuniqueness, there might be 
 energy solutions which are not limiting. Even though
 it is not excluded that (some of) these solutions may have 
 the same good regularity properties of the limiting 
 ones, they have to be forcedly excluded 
 from the long time analysis.
\eddos
To study the long time dynamics of ``limiting solutions'',
we also need to introduce a suitable 
extension (cf.~\cite[Sec.~2.2]{RSS})
of the concept of ``generalized
semiflow'' introduced by J.~Ball in the celebrated 
papers~\cite{Ba1,Ba2}.
\bede\label{defilimsem}
 We say that a family $\calS$ of maps from $[0,\infty)$ to
 a metric space $X$ is a\/ {\rm limiting semiflow} on the\/
 {\rm phase space} $X$ if the following properties 
 hold:
 \begin{enumerate}
 \item[{\rm (H1)}]
 For all $u_0\in X$ there exists at least one $u\in\calS$ such
 that $u(0)=u_0$\/ {\rm (existence property)}.
 \item[{\rm (H2)}]
 For all $u\in\calS$ and every $\tau\ge0$, the function $u^\tau$
 defined for $t\ge 0$ by $u^\tau(t):=u(t+\tau)$ still 
 belongs to $\calS$\/ {\rm (translation invariance)}.
 \item[{\rm (H4)}]
 For all sequence $(u_k)\subset \calS$ such that $u\zkk:=u_k(0)$ tends
 in $X$ to some $u_0$, there exist $u\in\calS$ such that $u(0)=u_0$
 and a nonrelabelled subsequence of $k$ such that, for all $t\ge 0$,
 it is $u_k(t)\to u(t)$ in $X$\/
 {\rm (upper semicontinuity w.r.t.~initial data)}.
 \end{enumerate}
\edde
\beos\label{noH3}
 In the above Definition, we kept Ball's notation. The property 
 missing here (obviously noted as (H3) in Ball's papers)
 essentially states that if we {\sl concatenate}\/ a couple of solutions 
 $u_1$, $u_2$, respectively defined on $[0,t_1]$ and on $[t_1,t_2]$
 and such that $u_1(t_1)=u_2(t_1)$, we still obtain
 a solution. Unfortunately, such a condition can be hardly proved
 for semiflows constructed by means of an approximation-limit 
 argument, and the main reason is that the values in $t_1$
 of the trajectories approximating the two solutions may be incompatible.
 From the point of view of asymptotics, it is shown 
 in \cite[Thm.~2.9]{RSS} (reported here as 
 Theorem~\ref{teoattrasemi} below) that the global attractor $\calA$ 
 still exists and is unique for limiting semiflows,
 under the natural conditions of dissipativity and eventual
 compactness. However, $\calA$ turns out to be just 
 {\sl quasi invariant}\/ (i.e., any $u_0\in \calA$ can be seen
 as the initial value of a {\sl complete orbit}\/ of the semiflow
 taking its values in $\calA$ for all $t\in\RR$) rather 
 than {\sl fully invariant}\/ (as it is in Ball's theory). 
 We refer to \cite[Subsec.~2.4]{RSS} for further remarks 
 on this point and to \cite{MV} for an alternative approach
 based on multivalued semiflows.
\eddos
The limiting solutions to \Pee~and (P$_0$) turn out to
fit the Definition above, at least
for a restricted class of potentials, including also
some {\sl singular}\/ cases:
\bede\label{defsepW}
 We say that the potential $W$ is\/ {\rm separating} if
 the following conditions are fulfilled. First, 
 \eqref{singW} holds. Second, for all $v\in \calW$ it is 
 $\max\{|v(x)|,x\in\barO\}<1$. Third, there exists
 an increasingly monotone function $\phi:[0,\infty)\to[0,1)$
 such that 
 \beeq{sepW}
   \|v\|_{C^0(\barO)}^2\le 
    \phi\big(\|v\|_{W^{1,6}(\Omega)}^2
             +\|\beta(v)\|^2\big)
     \quad\perogni v\in \calW.
 \end{equation}
\edde
\beos
 An easy refinement of the argument in 
 \cite[Prop.~2.10]{SSS1} shows that a sufficient
 condition for $W$ to be separating is that $W'$
 explodes sufficiently fast in proximity of $\pm1$. 
 Namely, in 3D, $W$ is separating provided that,
 for some $c>0$ (recall that $\beta=W'+\lambda\Id$),
 \beeq{quandoWsepa}
   \beta(r)\ge \frac c{(1-r)^3},\qquad
    -\beta(r)\ge \frac c{(1+r)^3},
 \end{equation}
 respectively in a left neighbourhood of $1$ and 
 in a right neighbourhood of $-1$. Actually, it is shown
 in \cite[Prop.~2.10]{SSS1} that, if $v\in\calW$
 (which entails that the 
 argument of $\phi$ in \eqref{sepW} is finite thanks to the 
 continuous embedding $H^2(\Omega)\subset W^{1,6}(\Omega)$)
 and \eqref{quandoWsepa} holds,
 then the maximum on $\barO$ of $v$ is strictly 
 smaller than $1$, in a way that (monotonically) depends on the 
 distance $\dW(v,0)$. Of course, we cannot exclude that
 a refinement of the argument in \cite[Prop.~2.10]{SSS1}
 might show that the class of separating potentials
 is in fact wider.
\eddos
The key point to show existence of the global attractor
is the following ``local regularization'' property of 
limiting solutions:
\bele\label{lemmalocreg}
 Let\/ \eqref{hpb}, \eqref{W}--\eqref{coercW}, and 
 \eqref{regof} hold. Let also $W$ either
 satisfy\/ \eqref{growW} for $p=\infty$ 
 or be\/ {\rm separating}. 
 Let $u$ be a limiting solution either to \Pee\ or 
 to\/ {\rm (P$_0$)}.
 Then, there exist $T_0$ depending 
 on $\calE(u_0)$, and $C_0,\delta>0$ 
 independent of $T_0$, $u_0$, and $\epsi$,
 such that, for all $T\ge T_0$ there holds the property
 \beeq{locsempreinW}
   \esiste~\tau=\tau(T)\in[0,3/2]:
    ~~\dW(u(t),0)\le C_0 \quad\perogni t\in[T+\tau,T+\tau+\delta].
 \end{equation}
\enle
Condition \eqref{locsempreinW} states that we are not able
to prove {\sl uniform in time}\/ regularization for limiting
solutions. Nevertheless, at least for a sequence of 
intervals whose length $\delta$ is uniformly controlled from 
below, the limiting solutions take values in a bounded 
ball of $\calW$, which is a relatively compact set of $\calV$
thanks to Lemma~\ref{lemmaWcompV}.
In this sense, \eqref{locsempreinW} can be thought as a 
``locally uniform'' regularization property.
We can now state the
\bete\label{teosemiflow}
 Let\/ \eqref{hpb}, \eqref{W}--\eqref{coercW}, and 
 \eqref{regof} hold. Let also $W$ either satisfy\/
 \eqref{growW} for $p=\infty$ or be\/ {\rm separating}. 
 Then, the\/ {\rm limiting solutions} to~\Pee\
 (respectively to\/ {\rm (P$_0$)}\/)
 constitute a {\rm limiting semiflow} $\calSe$
 (respectively $\calS_0$) on $\calVm$.
\ente
In the following statement \cite[Thm.~2.9]{RSS} (see
also \cite[Thm.~3.3]{Ba1}) we collect the definition
of global attractor for a limiting semiflow
and the basic tool to prove its existence.
\bete\label{teoattrasemi}
 Let $\calS$ be a limiting semiflow on the metric space $X$.
 Then, a compact set $\calA\subset X$ is the global attractor
 for $\calS$ if it is compact, quasi-invariant, and it
 attracts all bounded sets of $X$ w.r.t.~its metric. The 
 attractor $\calA$ exists if and only if $\calS$ satisfies
 the following properties:
 \begin{enumerate}
 \item[{\rm (A1)}]
 There exists a metric bounded set $B_0\subset X$ such that 
 any $u\in \calS$ eventually takes values in
 $B_0$\/ {\rm (point dissipativity)}.
 \item[{\rm (A2)}]
 For all $(u_k)\subset\calS$ such that $(u\zkk)$ is 
 bounded in $X$ (where $u\zkk:=u_k(0)$)
 and all $(t_k)$ such that
 $t_k\nearrow\infty$, there exist $u_\infty\in X$ and
 a (nonrelabelled) subsequence of $k$ such that 
 $u_k(t_k)\to u_\infty$ in $X$\/
 {\rm (asymptotic compactness)}.
 \end{enumerate}
 Finally, if $\calA$ exists, it is then unique.
\ente
%
%
Property \eqref{locsempreinW} and a careful use of the energy
equality are also the key tools to prove the 
\bete\label{teoattra}
 Let\/ \eqref{hpb}, \eqref{W}--\eqref{coercW}, and 
 \eqref{regof} hold. Let also $W$ either satisfy\/
 \eqref{growW} for $p=\infty$ or be\/ {\rm separating}. 
 Then, the semiflows $\calSe$, $\calS_0$ admit 
 compact global attractors $\calA_\epsi$, $\calA_0$ 
 in the sense of\/ {\rm Theorem~\ref{teoattrasemi}}. 
 Moreover, $\calA_\epsi$ is a metric bounded subset of 
 $\calWm$ if\/ \eqref{growW} holds with $p\le 6$.
\ente
%


\section{Proofs}
\label{secvisc}

In what follows, the symbols $c$, $c_i$, and $C_i$, 
with $i\ge0$, will denote positive constants, depending on the
data $b,W,f$ of the problem, but independent
of $\epsi,u_0$, of time, and of approximation
parameters (e.g., of $n$ in Problems~\Pnee,~\Pneen).
In particular, small letters $c$ and $c_i$ will
be used in the computations, and capital
letters $C_i$ in the resulting estimates.
Dependence on $m$ (cf.~\eqref{regou0})
is allowed. Moreover, the value of $c$ may
vary even inside a single line. The symbol
$c\OO$ will denote some embedding constants
depending only on $\Omega$. Finally,
$\cc$, $\cc_i$, $i\ge0$, will stand for positive
constants with additional dependences 
(e.g., on time or on $u_0$), specified time by time.

\vspace{2mm}

\noindent%
{\bf Proof of Theorem~\ref{teoesivisco}.}~~%
Let, for $n\in\NN$, $u$ be a solution either to 
Problem~\Pneen\ or to Problem~\Pnee\
introduced in the previous 
Section. The subscript $n$ is omitted in the
notation of $u$ just for brevity.
We now perform some a priori estimates, with the 
purpose of removing the $n$-approximation. 
Note that we can take advantage of the 
regularity \eqref{regoun} so that all the procedure
below makes sense rigorously.
Moreover, there holds
\beeq{propEn}
  \calE(v),\calEn(v)\ge\eta\|v\|_V^2-c
   \quad\perogni v\in V,~n\in\NN,
\end{equation}
where $\eta>0$ depends on $\lambda, c_W, f$ and 
is independent of $n$.
Testing \eqref{eqne} by $w$, 
\eqref{eqne2} by $u_t$, taking the sum, and using
\eqref{hpb}, we obtain the approximate 
energy equality
\beeq{energy}
  \ddt \calE_n(u)
   +\io \big(b(u)|\nabla w|^2\big)
   +\epsi\|u_t\|^2= 0,
\end{equation}
which holds at least
a.e.~in time. Hence, by \eqref{hpb},
$\calE_n$ is a Liapounov functional.
To get dissipativity, we also 
have to test \eqref{eqne2} by $u-u\OO$. 
By \eqref{medie1}, \eqref{consu} and the Young and 
Poincar\'e-Wirtinger inequality
(in the form
\beeq{PoWi}
  \|v-v\OO\|^2\le 
   c_\Omega\|\nabla v\|^2
   \quad \perogni v\in V
\end{equation}
and for some $c\OO>0$), we get
\bealo
  & \epsi\ddt\|u-u\OO\|^2
   +2\|\nabla u\|^2
   +2\io (W_n'(u)+f)(u-u\OO)
   =2\io w(u-u\OO)=2\io (w-w\OO)(u-u\OO)\\
 \label{conto10.2}
 & \mbox{}~~~~~~~~~~
   \le 2\|u-u\OO\|\|w-w\OO\|
   \le \|\nabla u\|^2+c_\Omega^2\|\nabla w\|^2.
\end{align}
Let us observe that,
by monotonicity of $r\mapsto \beta_n(r)=W_n'(r)+\lambda r$,
%
\beeq{WW'}
  W_n'(r)(r-r_0)
   \ge W_n(r)-W_n(r_0)-\lambda(r^2+r_0^2)
   \quad\perogni r,r_0\in \RR.
\end{equation}
Thus, noting that, thanks also to \eqref{regou0},
\beeq{WWn}
  W_n(u\OO)\le W(u\OO)
   \le c=c(W,m),
\end{equation}
and using \eqref{Wnewappr} and \eqref{PoWi},
the last term on the \lhs\ of 
\eqref{conto10.2} is estimated by
\bealo
  2\io(W_n'(u)+f)(u-u\OO)
   & \ge 2\io W_n(u)-2\lambda\|u\|^2
   +\io f(u-u\OO)-c\\
 \label{conto10.3}
   & \ge \io W_n(u)
   -\frac12\|\nabla u\|^2-c.
\end{align}
Thus, summing together \eqref{energy} and 
$2^{-1}\alpha c_\Omega^{-2}\times$\eqref{conto10.2},
and taking \eqref{hpb}, \eqref{conto10.3} into account,
we readily get, for some $\kappa,C_0>0$ with the 
same dependencies as the generic $c$,
\beeq{st11}
  \ddt \Big(\calEn(u)+\frac{\alpha\epsi}{2c\OO^2}\|u-u\OO\|^2\Big)
   +\kappa\big(\calEn(u)
   +\|\nabla w\|^2
   +\epsi\|u_t\|^2\big)
  \le C_0.
\end{equation}
By Gronwall's Lemma, this gives
\eqref{dissepsi} (cf.~Lemma~\ref{enerbuone}),
with $\calE_n$ in place of $\calE$.
Next, we test \eqref{eqne2} by $Bu$.
Using \eqref{coercW} together with
the H\"older and Young inequalities, we infer
\beeq{conto11}
  \epsi\ddt\|\nabla u\|^2
   +\|Bu\|^2
 \le\|f\|^2
  +(2\lambda+1)\|\nabla u\|^2
  +\|\nabla w\|^2.
\end{equation}
The combination of \eqref{energy} and \eqref{conto11}
then immediately gives \eqref{regou} (at the $n$-approximated
level).

To get \eqref{regow}, it remains to estimate the space 
averages of $W_n'(u)$ and $w\OO$. 
To do this, we can proceed by using an
argument devised in \cite{KNP} (see also \cite[Sec.~5]{BDS})
which is just sketched here. Namely, we first have to compute
%
%
%
\eqref{eqne2} times $\beta_n(u)-(\beta_n(u))\OO$.
By standard calculations, this gives
\beeq{contoknp1}
  \big\|\beta_n(u)-(\beta_n(u))\OO\big\|^2\le   
   c\big(1+\|\nabla u\|^2+\|\nabla w\|^2+\epsi^2\|u_t\|^2\big).
\end{equation}
Proceeding, e.g., as in \cite[(5.32)--(5.33)]{BDS},
we also get 
\beeq{contoknp2}
  \big|(\beta_n(u))\OO\big|^2\le 
   c\big\|\beta_n(u)-(\beta_n(u))\OO\big\|^2\|u-u\OO\|^2,
\end{equation}
where the constant $c$ depends in particular on $m$.
Then, the coupling of \eqref{contoknp1} and \eqref{contoknp2},
together with \eqref{regou} and a further comparison 
in \eqref{eqne2} (made in order to estimate $w\OO$),
readily give \eqref{regow} (note that all constants
in the procedure are independent of $n$). More 
precisely, integrating \eqref{energy},
\eqref{conto11} and \eqref{contoknp1}
from $t$ to $t+1$, for $t\ge\tau\ge 0$,
recalling \eqref{dissepsi},
and taking also \eqref{contoknp2} 
into account, it is not difficult to infer that,
for some monotone function $M:[0,\infty)\to [0,\infty)$
independent of $n$ and possibly new values of $C_0,\kappa$, 
it is 
\beeq{stunif1}
  \sup_{t\in(\tau,\infty)}
   \int_t^{t+1}\big(\|u(s)\|_{H^2(\Omega)}^2
         +\|w(s)\|_V^2
         +\|\beta_n(u(s))\|^2\big)\,\dis
    \le M(\calE(u_0))e^{-\kappa\tau}+C_0.
\end{equation}
To proceed, we now prove that, as at least 
a subsequence of $n$ goes to infinity,
the solutions to \Pnee\ (or \Pneen) pass to the
limit yielding (at least) one solution to~\Pee\ 
(or to (P$_0$)) which satisfies the
same bounds (and w.r.t.~the same constants). This is
standard and can be performed essentially as in \cite{BB} or
\cite{BDS}, so we shall give very few details.

Actually, the uniform bounds corresponding to
\eqref{regou} and \eqref{regow} 
entail that a for a (not relabelled) subsequence 
of $n$ there holds convergence to a proper 
limit function. 
Then, to show that this limit function solves
\eqref{eqne}--\eqref{eqne2} in the original form
and satisfies \eqref{iniz} one has to pass 
to the limit the nonlinear terms. In particular,
to identify the term depending on $W'$,
the G-convergence $W_n'\to W'$ and a standard
monotonicity argument are used (note that 
$W'$ is monotone up to a linear perturbation
and refer, e.g., to~\cite[Prop.~3.59, p.~361]{At} or
\cite[Prop.~1.1, p.~42]{Ba}).
Moreover, to treat the mobility term, one can
observe that by \eqref{hpb} $b(u)$ converges
strongly in $L^a(\Omega\times(0,T))$ for 
all $a\in[1,\infty)$ and weakly 
star in $L^\infty(\Omega\times(0,T))$ 
for all $T>0$. Thus, by the bound on $w$ in 
\eqref{stunif1}, also the product
$b(u)\nabla w$ (or, analogously, the 
operator $B_u$) passes to the proper limit.
Let us also notice that the dissipativity 
bounds \eqref{dissepsi} and \eqref{stunif1} 
are still valid in the limit with the same constant
$C_0$ and $\kappa$ thanks to semicontinuity 
property of norms w.r.t.~weak or weak star convergences.

Next, let us show (conditional) uniqueness. To start with,
note that if $u$ is an energy solution additionally
satisfying \eqref{regou2} for some $\tau\ge 0$, $T\ge\tau$,
then, by the first of \eqref{regou2}, the multiplication
operator
\beeq{multi}
  V\to V, \qquad v\mapsto b(u(t))v 
\end{equation}
is continuous for a.a.~$t\in (0,T)$.
Then, \eqref{eqne} can be rewritten, a.e.~in~$(\tau,T)$,
as the relation in $V'$ 
\beeq{eqneell}
  Bw=-\frac1{b(u)}u_t
   +\frac{b'(u)}{b(u)}\nabla u\cdot\nabla w.
\end{equation}
Moreover, evaluating the $H$-norm of the latter term
on the \rhs\ we have
\bealo
  \Big\|\frac{b'(u)}{b(u)}\nabla u\cdot\nabla w\Big\|
   & \le c \|\nabla u\|_{L^6(\Omega)}
       \|\nabla w\|_{L^3(\Omega)}
     \le c\|u\|_{H^2(\Omega)}\|\nabla w\|_{L^6(\Omega)}^{1/2}
       \|\nabla  w\|^{1/2}\\
 \label{conto21}
   & \le \cc \|w\|_{H^2(\Omega)}^{1/2} \|\nabla w\|^{1/2}
    \le \frac12\|Bw\|+ \cc\|w\|_V,
\end{align}
where the constants $\cc$ depend also on the
$L^\infty(\tau,T;H^2(\Omega))$-norm of $u$. Thus,
\beeq{addiregw}
  w\in L^2(\tau,T;H^2(\Omega)).
\end{equation}
Then, closely with the procedure 
in \cite[Proof of Thm.~2.2]{BB}, we consider a 
pair $u_1,u_2$ of solutions both satisfying~\eqref{regou2} 
(and consequently~\eqref{addiregw}) for some $\tau\ge 0$, 
$T\ge\tau$. We also assume that
$u_1(\tau)=u_2(\tau)$. Then, we set $u:=u_1-u_2$,
$w:=w_1-w_2$, compute 
$\big($\eqref{eqne}$_{1}-$\eqref{eqne}$_{2}\big)
\times \calN_{u_1} u$, and subtract 
$\big($\eqref{eqne2}$_{1}-$\eqref{eqne2}$_{2}\big)
\times u$ from the result.
Note that, actually, $u\in V_0$ a.e.~in
time. Setting $\zeta:=\calN_{u_1}u$, so that
$B_{u_1}\zeta=u$, we have 
\beeq{conto22}
  \frac\epsi2\ddt\|u\|^2
   +\|\nabla u\|^2
   +\big(\zeta,(B_{u_1}\zeta)_t\big)
  \le \lambda\|u\|^2
   -\io\big(b(u_1)-b(u_2)\big)\nabla w_2\cdot \nabla\zeta.
\end{equation}
Then, it is not difficult to see that 
\beeq{conto23}
  \big(\zeta,(B_{u_1}\zeta)_t\big)
   \ge \ddt\io\frac{b(u_1)}2|\nabla\zeta|^2
   -c\|u_{1,t}\|_{L^6(\Omega)}\|\nabla\zeta\|_{L^{12/5}(\Omega)}^2.
\end{equation}
Moreover, noting that 
there holds (cf.~\cite[(2.26)]{BB} or 
\cite[Lemma~1.2]{SchS})
\beeq{scheidschim}
  \|\zeta\|_{H^2(\Omega)}\le c\big(\|u\|
   +\|u_1\|_{H^2(\Omega)}^2\|\nabla \zeta\|\big),
\end{equation}
we get
\bealo
  \|\nabla\zeta\|_{L^{12/5}(\Omega)}^2
  & \le c\|\nabla\zeta\|^{3/2}\|\zeta\|_{H^2(\Omega)}^{1/2}
   \le c\|\nabla\zeta\|^{3/2}
  \big(\|u\|^{1/2}+\|u_1\|_{H^2(\Omega)}\|\nabla\zeta\|^{1/2}\big)\\
 \label{conto23.2}
   & \le \cc\|\nabla\zeta\|^2 + \cc_0\|u\|^2,
\end{align}
the constants $\cc,\cc_0$ depending also on $u_1$.
Finally, using \eqref{coerc},
the analogue of \eqref{scheidschim}, and 
standard embedding and interpolation inequalities,
we get
\bealo
  & -\io(b(u_1)-b(u_2))\nabla w_2\cdot \nabla\zeta
    \le c\|u\|_{L^3(\Omega)}\|\nabla w_2\|_{L^6(\Omega)}
          \|\nabla \zeta\|
    \le c\|\nabla u\|^{3/4}\|u\|_{V'}^{1/4}
     \|w_2\|_{H^2(\Omega)} \|\nabla \zeta\|\\
 \label{conto24}
 & \mbox{}~~~~~~~~~~
   \le c\|\nabla u\|^{3/4}
   \|w_2\|_{H^2(\Omega)} \|\nabla \zeta\|^{5/4}
    \le \frac14\|\nabla u\|^{2}
     +c\|w_2\|_{H^2(\Omega)}^{8/5}\|\nabla \zeta\|^2.
\end{align}
Thus, noting that, by interpolation,
\beeq{disucompa0}
  (\cc_0+\lambda)\|u\|^2
   \le \frac14\|\nabla u\|^2 + \cc\|\nabla\zeta\|^2,
\end{equation}
and since $w_2$ complies with \eqref{addiregw},
an application of the Gronwall Lemma in \eqref{conto22}
permits to conclude for uniqueness of regularizing solutions.

Finally, let us now assume $\epsi>0$ and that
\eqref{growW} holds with $p=6$
and let us prove \eqref{regou2}.
Of course, to be fully rigorous,
we should work on the solution to~\Pnee\ and then 
pass to the limit, but, for brevity and since 
everything is standard, we assume here 
directly that $u$ solves~\Pee. Then, 
testing \eqref{eqne2} by $Bu_t$, we obtain
\bealo
  & \ddt \big(\|Bu\|^2+2(f,Bu)\big)
   +2\epsi\|\nabla u_t\|^2
  \le 2(w-W_n'(u),Bu_t)\\
 \label{conto12}
  & \mbox{}~~~~~ \le  \epsi\|\nabla u_t\|^2
   +\frac1\epsi\|\nabla w+W_n''(u)\nabla u\|^2.
\end{align}
To estimate the latter term, we use \eqref{growW}
(which holds uniformly in $n$), \eqref{regou}, standard 
interpolation inequalities, and the continuous
embeddings $V\subset L^{6}(\Omega)$ and 
$H^{5/4}(\Omega)\subset L^{12}(\Omega)$.
Namely, we have
\bealo
  \|W_n''(u)\nabla u\|^2
   & \le \|\nabla u\|_{L^6(\Omega)}^2\big(1+\|u\|_{L^{12}(\Omega)}^8\big)\\[1mm]
  \label{conto12bis}
   & \le \|u\|_{H^2(\Omega)}^2\big(1+\|u\|_{V}^6\|u\|_{H^2(\Omega)}^2\big).
\end{align}
Thus, \eqref{conto12} gives, for all $t\ge\tau\ge0$,
\beeq{conto13}
  \ddt \big(\|Bu\|^2+2(f,Bu)\big)
   +\epsi\|\nabla u_t\|^2
  \le c\ee\big(M(\calE(u_0))e^{-\kappa\tau}+C_0\big)
   \big(1+\|Bu\|^4\big),
\end{equation}
where $c\ee>0$ depends on $\epsi$ and it explodes 
as $\epsi\searrow0$.
Thus, recalling \eqref{stunif1} and using the 
standard Gronwall Lemma if it is $u_0\in\calW_m$ and 
the {\sl uniform}\/ Gronwall Lemma \cite[Lemma~I.1.1]{Te},
if it is just $u_0\in\calV_m$, we get the regularity
of $u$ in \eqref{regou2} respectively for $\tau\ge 0$ 
and $\tau>0$. To conclude, we notice that,
being $p=6$ in~\eqref{growW},
$\|\beta(u)\|\le c(1+\|u\|_{L^8(\Omega)}^4$),
so that also the third of \eqref{regou2} holds.
%
%
The proof of Theorem~\ref{teoesivisco} is thus
complete.\dimbox

\vspace{2mm}

\noindent%
{\bf Proof of Lemma~\ref{enerbuone}.}~~%
The energy equality can be proved as for the approximated
problem (i.e., testing \eqref{eqne} by $w$, \eqref{eqne2}
by $u_t$, and taking the difference). The key point is
that, for $\epsi>0$, the regularities 
\eqref{regou}--\eqref{regow} are sufficient to apply
the chain rule formula~\cite[Lemma~3.3, p.~73]{Br};
actually, all the terms in \eqref{eqne2} lie in $\LDH$
as well as the test function $u_t$. As we consider, instead, 
energy solutions to (P$_0$), $u_t$ is only in $\LDVp$
and the terms in \eqref{eqne2} do not lie, each one
separately, in $\LDV$. Nevertheless, since both $w$ 
and the sum $Bu+W'(u)+f$ lie in in $\LDV$, one can 
use, e.g., \cite[Lemma~4.1]{RS} (note that $W$ satisfies
the growth assumption \cite[(4.23)]{RS} due to
\eqref{Wnew}) and still conclude for \eqref{energyeq}.
Observe also that, as a byproduct, the function
$t\mapsto \calE(u(t))$ is absolutely continuous on
$[0,T]$. Finally, by the same type of considerations,
also the procedure used to get~\eqref{dissepsi}
(cf.~the computation leading
to~\eqref{st11}) can be justified.\dimbox 
%
%
%

\vspace{2mm}

\noindent%
{\bf Proof of Lemma~\ref{lemmalocreg}.}~~%
Let $u$ be a limiting solution 
either to Problem~\Pee\ or to (P$_0$), and 
let $u_j$, $n_j$ and $\epsi_j$ be as in
Definition~\ref{defisollimi}. Being 
$u_j\in \calY_T$ for all $j\in\NN$ and $T>0$, 
we can test \eqref{eqne2j} by $Bu\jjt$. 
We then get (cf.~\eqref{conto12}) 
\beeq{conto30}
  \ddt \big(\|Bu_j\|^2+2(f,Bu_j)\big)
   +2\epsi_j\|\nabla u\jjt\|^2
  \le 2(w_j-\Wnj'(u_j),Bu\jjt).
\end{equation}
Next, we test \eqref{eqnej} by $w\jjt$ and subtract
from the result the expression obtained by 
differentiating in time \eqref{eqne2j} and testing
it by $u\jjt$.
Using \eqref{hpb} and 
\eqref{coercW} (which still holds for $\Wnj$),
we get
\beeq{conto31}
  \ddt \Big(\epsi_j\|u\jjt\|^2
   +\io \big(b(u_j)|\nabla w_j|^2\big)\Big)
   +2\|\nabla u\jjt\|^2
  \le 2\lambda\|u\jjt\|^2
   +c\io\big(|u\jjt||\nabla w_j|^2\big).
\end{equation}
Next, we compute
\bealo
  \ddt \big( \|u_j\|^2 + \|\beta_{n_j}(u_j)\|^2 \big)
   & = 2 (u_j,u\jjt)
    + \io \big(\beta_{n_j}'(u_j)\beta_{n_j}(u_j)u\jjt\big)\\
 \label{conto32}
   &\le  \|u_j\|^2+\|u\jjt\|^2
    + \io \big(\beta_{n_j}'(u_j)\beta_{n_j}(u_j)u\jjt\big).
\end{align}
Then, we estimate the terms on the 
\rhs s of \eqref{conto30} and \eqref{conto31}:
\beal{conto33a}
  2(w_j,Bu\jjt) 
   & \le 4 \|\nabla w_j\|^2 + \frac14 \|\nabla u\jjt\|^2,\\
 \no
  -2(\Wnj'(u_j),Bu\jjt) 
   & =2\lambda(u_j,Bu\jjt)-2(\beta_{n_j}(u_j),Bu\jjt)\\
 \label{conto33b}
   & \le \frac14\|\nabla u\jjt\|^2 + c\|\nabla u_j\|^2
      - 2 \io\big(\beta_{n_j}'(u_j)\nabla u_j\cdot \nabla u\jjt\big),\\
 \label{conto33c}
  \io\big(|u\jjt||\nabla w_j|^2\big)
   &  \le c\|u\jjt\|_{L^6(\Omega)}\|\nabla w_j\|^{3/2}
     \big(\|\nabla w_j\|^{1/2}+\|Bw_j\|^{1/2}\big).
\end{align}
Using the analogue of
\eqref{scheidschim}, one gets from \eqref{conto33c} that
\beeq{conto34}
  \io\big(|u\jjt||\nabla w_j|^2\big)
   \le \frac14\|u\jjt\|_V^2
    +c\|\nabla w_j\|^4+c\|\nabla w_j\|^6
    +c\|u_j\|^6_{H^2(\Omega)}.
\end{equation}
Next, we estimate the last terms on the \rhs s
of \eqref{conto32} and \eqref{conto33b} this way:
\bealo
  & 2\io \beta_{n_j}'(u_j)\big(\beta_{n_j}(u_j)u\jjt
   -\nabla u_j\cdot \nabla u\jjt\big)\\
 \label{conto35}
  & \mbox{}~~~~~~~~~~
  \le \frac14\|u\jjt\|_V^2
   +c\|u_j\|_{H^2(\Omega)}^6
   +c\|\beta_{n_j}(u_j)\|^6
   +c\|\beta_{n_j}'(u_j)\|_{L^3(\Omega)}^3.
\end{align}
Summing now \eqref{conto30}, \eqref{conto31} and  
\eqref{conto32}, using on $H^2(\Omega)$ the equivalent norm
$(\|\cdot\|^2+\|B\cdot\|^2)^{1/2}$, and owing to 
\eqref{conto33a}--\eqref{conto35}, we get
\bealo
  & \ddt \Big(\|u_j\|^2_{H^2(\Omega)}
   +2(f,Bu_j)
   +\epsi_j\|u\jjt\|^2 
   +\io \big(b(u_j)|\nabla w_j|^2\big)
   +\|\beta_{n_j}(u_j)\|^2\Big)
   +(1+2\epsi_j)\|\nabla u\jjt\|^2\\
 \label{conto36}
  & \mbox{}~~~~~~~~~~ 
  \le c_1 \big(1 
   +\|u\jjt\|^2
   +\|\nabla w_j\|^6
   +\|u_j\|_{H^2(\Omega)}^6
   +\|\beta_{n_j}(u_j)\|^6
   +\|\beta_{n_j}'(u_j)\|_{L^3(\Omega)}^3 \big).
\end{align}
Again, by interpolation, and recalling \eqref{hpb},
we have
\beeq{disucompa}
  c_1\|u\jjt\|^2
   \le \frac12\|\nabla u\jjt\|^2 + c\|u\jjt\|_{V'}^2
   \le \frac 12\|\nabla u\jjt\|^2 + c\|\nabla w_j\|^2.
\end{equation}
To proceed, we start considering the simpler 
case when $I=\RR$ (i.e., \eqref{growW} holds with
$p=\infty$). Then, noting as 
$\upsilon_n:=(\Id+n^{-1}\beta)^{-1}$
the {\sl resolvent}\/ of $\beta$ of index $n^{-1}$, 
and recalling that, for all $n$, $\upsilon_n$ is 
1-Lipschitz and satisfies
$\beta_n=\beta\circ \upsilon_n$
and $\upsilon_{n_j}(0)=0$ (due to \eqref{W}),
it is not difficult to realize that
\beeq{conto37}
  \|\beta_{n_j}'(u_j)\|_{L^3(\Omega)}^3
   \le \|\beta'(\upsilon_{n_j}(u_j))\|_{L^3(\Omega)}^3 
   \le \gamma\big(\|u_j\|^2_{L^{\infty}(\Omega)}\big),
 \end{equation}
where we have set, for $s\ge0$, 
\beeq{defigamma}
   \gamma_0(s):=|\Omega|\max\big\{\beta'(r)^3+\beta'(-r)^3,
     r\in [0,s]\big\}, \qquad
   \gamma(r):=\gamma_0\big(s^{1/2}\big),  
\end{equation}
so that $\gamma:[0,\infty)\to [0,\infty)$ is 
monotone. By continuity of the embedding $H^2(\Omega)\subset
L^\infty(\Omega)$, 
it follows that, for all $c_2\ge 0$,
\bealo
  & \ddt \Big(c_2+\|u_j\|^2_{H^2(\Omega)}
   +2(f,Bu_j)
   +\epsi_j\|u\jjt\|^2 
   +\io \big(b(u_j)|\nabla w_j|^2\big)
   +\|\beta_{n_j}(u_j)\|^2\Big)
   +\Big(\frac12+2\epsi_j\Big)\|\nabla u\jjt\|^2\\
 \label{conto38}
  & \mbox{}~~~~~~~~~~ 
  \le c_1\Big(1+\|\nabla w_j\|^6
   +\|u_j\|^6_{H^2(\Omega)}
   +\|\beta_{n_j}(u_j)\|^6
   +\gamma\big(c\OO\|u_j\|_{H^2(\Omega)}^2\big)\Big).
\end{align}
Then, noting as $y_j$ the function whose time 
derivative appears on the \lhs\ and 
choosing then $c_2>0$ so that $y_j\ge 0$
(note this can be done independently 
of the initial datum), the relation above
can be interpreted in the form
\beeq{conto39}
  y_j'(t)\le \psi(y_j(t)),
\end{equation}
where $\psi:[0,\infty)\to[0,\infty)$ is a suitable
monotone function depending on $\gamma$, but 
independent of $j$. Now, let us observe that $u_j$ 
satisfies the analogue of \eqref{stunif1}, namely
\beeq{stunif1new}
  \sup_{t\in(T_0,\infty)}
   \int_t^{t+1}\big(\|u_j(s)\|_{H^2(\Omega)}^2
         +\|w_j(s)\|_V^2
         +\|\beta_{n_j}(u_j(s))\|^2\big)\,\dis
    \le M\big(\calE(u_0)+\sigma_j\big)e^{-\kappa T_0}+C_0,
\end{equation}
where $\sigma_j$ is as in \eqref{conve4} and $T_0>0$.
Actually, the above surely holds if $u_j$ is a solution to 
some~\Pnee\ (or~\Pneen). Here, although 
the choice of $u_j$ is slightly more 
general (cf.~Remark~\ref{nonproprioPn}), it is 
easy to see that \eqref{stunif1new} is still satisfied.
Then, taking $j_0$ so large that $\sigma_j\le \calE(u_0)$
for all $j\ge j_0$, it is clear that also $T_0$ can be 
chosen so that the \rhs\ above is $\le 2C_0$.
Thus, we can directly assume $T_0=1$ for simplicity 
of notation and without loss of generality.
In what follows, we note as $(j)_i$
a subsequence of $(j)_{j\ge j_0}=:(j)_0$ obtained
by successive extractions. Namely, for all $i$,
$(j)_i$ is a subsequence of $(j)_{i-1}$.
The indexes belonging to $(j)_i$ will be still indicated
just by $j$, for notational simplicity.

Relation \eqref{stunif1new} readily implies
that there exists $C_1>0$ 
such that, for all $m\in \NN$ (recall we assumed
$T_0=1$) and all $j\in(j)_0$
there exists $t\mmj\in [m,m+1/2]$ such that 
$y_j(t\mmj)\le C_1$. Then, defining
\beeq{defipsi}
  \Psi(s):=\int_{C_1}^s\frac{\dir}{\psi(r)}
\end{equation}
and solving \eqref{conto39}, it is clear that
\beeq{conto40}
  \Psi(y_j(t))=\Psi(y_j(t\mmj))
   +(t-t\mmj)
  \le t-t\mmj,
\end{equation}
at least for all $t\ge t\mmj$ such that
the relation above makes sense. This implies in
particular that there exist $\delta\in(0,1/4]$ and
$C_2>0$, both independent of $\epsi$,
$m$, $j$ and $u_0$ and such that
\beeq{conto41}
  \|u_j\|_{L^\infty(t\mmj,t\mmj+2\delta;H^2(\Omega))}^2
   +\|\beta_{n_j}(u_j)\|_{L^\infty(t\mmj,t\mmj+2\delta;H)}^2 
   \le \Psi^{-1}(\delta)=C_2<\infty,
\end{equation}
which holds for all $m\in\NN$ and $j\in(j)_0$.
Since the sequence $j\mapsto t_{1,j}$, $j\in(j)_0$ ranges in 
the compact interval $[1,3/2]$, it is clear that 
we can extract a subsequence $(j)_1$ and find a point 
$t_1\in [1,3/2]$, such that
\beeq{conto42}
  \|u_j\|_{L^\infty(t_i,t_i+\delta;H^2(\Omega))}^2
   +\|\beta_{n_j}(u_j)\|_{L^\infty(t_i,t_i+\delta;H)}^2 \le C_2
\end{equation}
holds for $i=1$ and for all $j\in(j)_1$.
Proceeding by induction,
for all $N\in\NN$ we then find $(j)_N$ such that 
\eqref{conto42}, where $t_i$ is some point in 
$[i,i+1/2]$, holds for all $i\le N$ and $j\in (j)_N$. 
At the end, we can thus extract a diagonal 
subsequence $(j)_\infty$, which gives \eqref{conto42} 
for all $j\in (j)_\infty$ and all $i\in\NN$.
Since $(j)_\infty$ is a subsequence of $(j)_0$, taking the 
$\liminf$ for $j\nearrow\infty$, $j\in(j)_\infty$,
of \eqref{conto42}, and recalling that,
by \eqref{conve3}, $\beta_{n_j}(u_j)$
tends to $\beta(u)$ weakly in $\LDH$ 
for all $T>0$, we finally get that $u$
satisfies the {\sl locally uniform}\/ regularization
estimate 
\beeq{conto43}
  \|u\|_{L^\infty(t_i,t_i+\delta;H^2(\Omega))}^2
   +\|\beta(u)\|_{L^\infty(t_i,t_i+\delta;H)}^2 \le C_2,
\end{equation}
with $t_i$ as above. This can be also rewritten as 
\beeq{locunif}
  \dW^2(u(t),0) \le C_2
   \qquext{for a.e.~}\, t\in\bigcup_{i=1}^\infty [t_i,t_i+\delta]
\end{equation}
and clearly entails the validity of \eqref{locsempreinW}
in case $W$ satisfies \eqref{growW} with $p=\infty$.

To conclude the proof, we have to face the case when 
$W$ is a {\sl separating potential}\/ 
(cf.~Definition~\ref{defsepW}). Then, the procedure
does not change till \eqref{disucompa}. 
After that point, the last inequality 
in \eqref{conto37} is replaced by
(notice that now $\upsilon_{n_j}$ takes
values into $(-1,1)$)
%
%
%
\bealo
  \|\beta_{n_j}'(\upsilon_{n_j}(u_j))\|_{L^3(\Omega)}^3
   & \le  \gamma\big(\|\upsilon_{n_j}(u_j)\|_{L^\infty(\Omega)}^2\big)
    \le   \gamma\Big(\phi\big(\|\upsilon_{n_j}(u_j)\|^2_{W^{1,6}(\Omega)}
       +\|\beta(\upsilon_{n_j}(u_j))\|^2\big)\Big)\\
 \label{conto45}
  & \le \zeta\big(\|u_j\|^2_{W^{1,6}(\Omega)}
       +\|\beta_{n_j}(u_j)\|^2\big)
     \le \zeta\big(c\OO\|u_j\|^2_{H^2(\Omega)}
       +\|\beta_{n_j}(u_j)\|^2\big),
\end{align}
where $\gamma$ is as in \eqref{defigamma},
$\phi$ as in Definition~\ref{defsepW}, 
and we have set $\zeta:[0,\infty)\to[0,\infty)$ as 
$\zeta:=\gamma\circ\phi$. Actually, still $\zeta$ is a
monotone function. Note that in \eqref{conto45}
we also used the 1-Lipschitz continuity of
$y_{n_j}$, Sobolev's embeddings,
and that $\upsilon_{n_j}(0)=0$ for all $j$.
At this point, one gets an expression similar
to \eqref{conto38}, where the last 
term $\gamma(c\OO\|u_n\|_{H^2(\Omega)}^2)$
is suitably replaced by the
\rhs\ of \eqref{conto45} (possibly up to a modification
of the expression of $\zeta$).
From this point on, the proof goes through with no 
further change.\dimbox

\vspace{2mm}

\noindent%
{\bf Proof of Theorem~\ref{teosemiflow}.}~~%
Property (H1) is an easy consequence of the Proof
of Theorem~\ref{teoesivisco}. In particular,
\eqref{conve1}--\eqref{conve2}, 
\eqref{conve3} and \eqref{conve4} follow from
\eqref{st11}, \eqref{conto11}, 
\eqref{contoknp1}--\eqref{contoknp2}, 
and \eqref{uzn}--\eqref{Wuzn}.
Property~(H2) is immediate. To show~(H4), let us 
first extract a nonrelabelled subsequence of $k$ such that,
for some function $u$, it is
\beeq{convinVT}
  u_k\to u \quext{weakly in }\,\calV_T \quad\perogni T>0.
\end{equation}
To do this, let us take $R>0$ such that $\dV(u\zkk,0)\le R$ 
for all $k$. Then, since any of the $(u_k)$ satisfies the 
analogue of \eqref{dissepsi}, setting $(k)_0:=(k)_{k\in\NN}$,
for all $N\ge 1$ we can extract a subsequence $(k)_N$ of
$(k)_{N-1}$ such that, for some function $u$, 
$u_k$ tends to $u$ at least weakly in $\calV_N$ 
as $k\in(k)_N$ goes to infinity. Thus, taking a diagonal
subsequence $(k)_\infty$,
it is clear that \eqref{convinVT} holds.
From this point on, we shall work on this subsequence.
Proceeding similarly with the Proof of 
Theorem~\ref{teoesivisco} (i.e., passing to the limit
in \eqref{eqne}--\eqref{eqne2}),
one sees immediately that
$u$ is an energy solution either to \Pee\ or to (P$_0$) 
and in particular it satisfies $u(0)=u_0$.
Actually, for all $T>0$ and $k\in(k)_\infty$,
it is, by the analogue of~\eqref{stunif1}, 
\beeq{nuovaW'}
  \|W'(u_k)\|_{L^2(0,T;H)}\le \cc
\end{equation}
(where $\cc$ might depend on $T$ and $R$) 
so that, by the {\sl strong}\/ convergence
$u_k\to u$ in $\LDH$, following from \eqref{convinVT},
and the usual monotonicity argument,
one has (without extracting any other subsequence
and for all $T>0$)
\beeq{nuovaW'2}
  W'(u_k)\to W'(u)\quext{weakly in }\,\LDH.
\end{equation}
Let us now show that $u$ is limiting, which is a bit 
more difficult. Since $u_k$ 
is limiting, for all $k\in (k)_\infty$
there exist an increasing sequence $j\mapsto n\jk$ 
(if $\epsi=0$, also a decreasing sequence
$j\mapsto\epsi\jk$, otherwise we intend that
$\epsi\jk\equiv\epsi$) and a sequence of 
functions $j\mapsto u\jk$, where $u\jk\in \calY_T$ 
for all $T>0$, solving, in the usual sense,
\beal{eqnek}
  & u\jkt+B_{u\jk} w\jk=0,\\
 \label{eqne2k}
  & w\jk=\epsi\jk u\jkt
      +Bu\jk+W_{n\jk}'(u\jk)+f.
\end{align}
Moreover, as $j\nearrow\infty$, $u\jk$ tends to $u_k$ 
in the sense specified in Definition~\ref{defisollimi}.
Then, it is clear that 
there exists $\cc>0$ depending on $R$ and $T$ 
but independent of $j$ and $k$, such that 
\beeq{unifnk} 
  \|u\jk\|_{\HUVp}
   +\epsi^{1/2}\|u\jk\|_{\HUH} 
   +\|u\jk\|_{\LIV} 
   +\|u\jk\|_{\LDHD} 
  \le \cc \qquad\perogni j,\,k.
\end{equation}
Next, for each $k\in(k)_\infty$
we can choose an index $j_k$ such that
the sequences $k\mapsto j_k$,
$k\mapsto n_k^{j_k}$ are strictly increasing
(and, if $\epsi=0$, $k\mapsto \epsi_k^{j_k}$ is strictly
decreasing) and 
\beeq{nkn}
  \|u\kjk-u_k\|_{C^0([0,k];H)}\le 1/k,
\end{equation}
so that, with no further extraction of subsequences
(the limit is already identified),
$k\mapsto u\kjk$ tends to the above constructed 
$u$ strongly in $C^0([0,T];H)$ and weakly
in $\calV_T$ for all $T>0$. This shows that
\eqref{conve1} and \eqref{conve2}, intended as
$j_k\nearrow\infty$, hold for the limit function 
$u$. Moreover, being for all $k$ (cf.~\eqref{conve4})
\beeq{provv1}
  \calE_{n\kjk}(u\kjk(0))\le \calE(u\zkk)+\sigma\kjk,
\end{equation}
it is clear that one can also take $(j_k)$ in such a
way that $j\mapsto\sigma\kjk$ is decreasing and tends to $0$, 
which shows \eqref{conve4} to hold for $u$ 
since $\calE(u\zkk)$ tends to $\calE(u_0)$ with $k$
by the hypothesis of convergence $u\zkk\to u_0$ in $\calV$
and thanks to Lemma~\ref{lemmaWcompV}. 
Next, noticing that, again,
\beeq{provv2}
  u\kjk\to u\quext{strongly in }\,\LDH,
   \qquad \big\|W'_{n\kjk}(u\kjk)\big\|_{\LDH}\le \cc,
\end{equation}
for all $T>0$, and using that $n\kjk\nearrow\infty$ with 
$k$ so that $W_{n\kjk}'$ G-converges to $W'$, 
one readily gets \eqref{conve3} for the limit $u$
by the usual monotonicity argument and still
for the whole sequence $(k)_\infty$.
Thus, $u$ is limiting.
%
%

The proof of (H4), however, is not yet concluded since
we still have to check that, choosing an
arbitrary $t\ge 0$, 
$u_k(t)$ tends to $u(t)$ {\sl strongly}\/
in $\calV$, which is {\sl not}\/ a consequence of 
the ``weak'' convergence in $\calV_T$
holding for all $T\ge 0$.
Actually, by the uniform bound on $u_k$ corresponding
to \eqref{unifnk}, it is clear that
$u_k\to u$ in $C_w(0,T;V)$ so that, for all $t\ge 0$,
one can only deduce that $u_k(t)$ tends to $u(t)$
{\sl weakly}\/ in $V$.

Thus, let us pick $T$ larger than the chosen $t$ 
and so large (in a way that only depends on $R$) that, 
by \eqref{locsempreinW}, for all $k$ in our subsequence
there exists $\tau_k\in [0,3/2]$ 
such that $\dW(u_k(s),0)\le C_0$ for all 
$s\in [T+\tau_k,T+\tau_k+\delta]$, where we remark
once more that $\delta$ and $C_0$ are 
independent of $k$ and $R$.

Now, the weak convergence in $\calV_T$
and \eqref{coerc} guarantee that 
\beeq{convwkw}
  w_k\to w \quext{weakly in }\,\LDV.
\end{equation}
%
%
Moreover, being $(\tau_k)\subset [0,3/2]$, 
which is a compact set, there exist
$S\in[T,T+3/2+\delta]$ and a subsequence $(k)_*$
of $(k)_\infty$ such that, at least for sufficiently
large $k\in(k)_*$, it is $\dW(u_k(S),0)\le C_0$. This, 
by Lemma~\ref{lemmaWcompV}, entails
that $u_k(S)$ tends to $u$ in $\calV$
so that, in particular, $\calE(u_k(S))$
tends to $\calE(u(S))$ at least as $k\in(k)_*$
goes to $\infty$. Next, 
writing the energy equality
\eqref{energyeq} for $u_k$ on the interval $(0,S)$
gives (possibly for $\epsi=0$)
\beeq{enuk}
  \calE(u_k(S))-\calE(u\zkk)
   =-\int_0^S\io b(u_k)|\nabla w_k|^2
   -\int_0^S\io \epsi|u_{k,t}|^2.
\end{equation}
Thus, taking the limit $k\nearrow\infty$
in $(k)_*$, noting that the \lhs\ converges to 
$\calE(u(S))-\calE(u_0)$, and using
the energy equality for $u$ we get 
(still possibly for $\epsi=0$)
\beeq{enuinfty}
  \lim_{k\nearrow\infty}\bigg(\int_0^S\io b(u_k)|\nabla w_k|^2
   +\int_0^S\io \epsi|u_{k,t}|^2\bigg)
  = \int_0^S\io b(u)|\nabla w|^2
   +\int_0^S\io \epsi|u_t|^2,
\end{equation}
which readily entails (in case $\epsi>0$, also thanks 
to the {\sl weak}\/ convergence
$u_{k,t}\to u_t$ in $L^2(0,S;H)$ following from
\eqref{convinVT})
\beeq{enuinfty2}
  \limsup_{k\nearrow\infty}\int_0^S\io b(u_k)|\nabla w_k|^2
   \le \int_0^S\io b(u)|\nabla w|^2.
\end{equation}
Let us now notice that, by \eqref{hpb},
\bealo
  \alpha\int_0^S\|\nabla w_k-\nabla w\|^2
   & \le \int_0^S\io\big(b(u_k)|\nabla w_k-\nabla w|^2\big)\\
 \label{enuinfty3}
  & =\int_0^S\io b(u_k)|\nabla w_k|^2
    +\int_0^S\io b(u_k)|\nabla w|^2
    -2\int_0^S\io b(u_k)\nabla w_k\cdot\nabla w.
\end{align}
Here, the latter two terms on the \rhs\
converge to the expected limits since 
$|\nabla w|^2\in L^1(\Omega\times(0,S))$
and there hold the convergences $b(u_k)\to b(u)$
(weakly star in $L^\infty(\Omega\times(0,S))$ and 
strongly in $L^a(\Omega\times(0,S))$ 
for all $a\in(1,\infty)$) 
and $\nabla w_k\to \nabla w$ 
(weakly in $L^2(0,S;H)$), which in 
particular entail 
\beeq{enuinfty2.2}
  b(u_k)\nabla w_k\to b(u)\nabla w \quext{{\sl weakly}\/ in }\,L^2(0,S;H).
\end{equation}
Thus, taking the $\limsup$ in \eqref{enuinfty3} 
and using \eqref{enuinfty2}, we obtain
\beeq{enuinfty4}
  w_k\to w \quext{{\sl strongly}\/ in }\,L^2(0,S;V),
\end{equation}
which, being $S\ge T$, implies in particular
\beeq{enuinfty4.2}
  w_k\to w \quext{{\sl strongly}\/ in }\,L^2(0,T;V).
\end{equation}
Notice that, a priori, the latter convergence holds 
only for the subsequence $(k)_*$ but, 
in fact, being the limit already identified, it is 
valid for the whole sequence $(k)_\infty$.
Thus, we can now come back to 
\eqref{enuk}, which we rewrite with $t$ in place 
of $S$. Using \eqref{enuinfty4.2} (and also \eqref{convinVT}
if $\epsi>0$), and recalling that $t\le T$,
we then get
\beeq{enuk2}
  \limsup_{k\nearrow\infty}\calE(u_k(t))
   \le \calE(u_0)
   -\int_0^t\io b(u)|\nabla w|^2
   -\int_0^t\io \epsi|u_{t}|^2,
\end{equation}
so that, by comparison in the limit energy equality
and thanks to Lemma~\ref{lemmaWcompV},
we finally get
\beeq{enuk3}
  \lim_{k\nearrow\infty}\calE(u_k(t))=\calE(u(t)),
\end{equation}
which implies that $u_k(t)\to u(t)$ 
in~$\calV$ and concludes the proof of~(H4) and
of the Theorem.\dimbox
\beos\label{perchecomp}
 The main reason which forced us to use the complicated
 ``local compactness'' argument is the presence of
 the nonconstant mobility $b(\cdot)$. Actually, 
 for constant $b$, once \eqref{convwkw} is known,
 one can immediately pass to \eqref{enuk2} and get 
 directly the strong convergence 
 $\calE(u_k(t))\to \calE(u(t))$. Instead, for 
 nonconstant $b$, without the 
 help of \eqref{enuinfty4.2} it is 
 not clear whether the semicontinuity property
 \beeq{enuinfty5}
  \int_0^t\io b(u)|\nabla w|^2
   \le\liminf_{k\nearrow\infty}\int_0^t\io b(u_k)|\nabla w_k|^2
 \end{equation}
 (which is necessary to prove \eqref{enuk2}) holds.
 Actually, at this stage, the integrand on the
 \rhs\ is only bounded in $L^1(\Omega)$ and  
 not even known to converge pointwise.
\eddos

\vspace{2mm}

\noindent%
{\bf Proof of Theorem~\ref{teoattra}.}~~%
Condition (A1) of Theorem~\ref{teoattrasemi} is an
immediate consequence of \eqref{dissepsi}. Let us then
show (A2). With the notation of (A2), let us first
point out that, being $(u\zkk)$ a bounded set in $\calV_m$,
by \eqref{locsempreinW} there exist $\delta\in[0,1/4]$ and 
$C_0>0$ such that for all (sufficiently large, depending on
the ``radius'' of $(u\zkk)$) $k\in\NN$ there 
exists $\tau_k\in[0,3/2]$ with
\beeq{attr1}
  \dW(u_k(t),0) \le C_0 \quad
   \perogni t\in [t_k-2+\tau_k,t_k-2+\tau_k+\delta].
\end{equation}
In particular, we can extract a subsequence of $(k)$,
not relabelled, such that $\tau_k\to \tau\in[0,3/2]$.
Setting then $v_k(s):=u_k(t_k-2+\tau+\delta/2+s)$,
and eventually being $|\tau_k-\tau|\le \delta/2$,
we then have that, still up to the extraction of 
a further subsequence, $v_k(0)$ tends to some 
$v_0$ 
in $\calV$. Moreover, by (H2),
$(v_k)\subset \calSe$ (possibly for $\epsi=0$), i.e.,
it is a limiting solution, and, clearly, 
\beeq{attr2}
  u_k(t_k)=v_k(2-\tau-\delta/2).
\end{equation}
Thus, the same argument used to show (H4) in the Proof
of Theorem~\ref{teosemiflow} permits to say that 
a subsequence of $v_k(2-\tau-\delta/2)$ admits a proper
limit (which coincides, by the way, with an element
of the semiflow evaluated at the time $2-\tau-\delta/2$)
in the metric topology of $\calV$.
This gives~(A2). 

Finally, if \eqref{growW} holds with $p\le 6$, it 
is clear from the uniform regularization property
\eqref{regou2} that $\calA\ee$ is bounded in $\calW_m$.
This concludes the proof of the Theorem.\dimbox


\section{Entropy solutions}
\label{secentropy}

In this Section, we show that, if $\epsi=0$,
\eqref{growW} holds with $p\in(2,6)$,
and, in place of \eqref{Wnew}--\eqref{coercW}, we have
\beeq{newW}
  W''(r)\ge\eta|r|^{p-2}-\lambda\qquad\perogni r\in I=\RR,
\end{equation}
where $\eta>0$ and $p\in(2,6)$ is the same exponent as in \eqref{growW},
then there exist weaker solutions to the analogue of 
Problem~(P$_0$), corresponding to
the choice of an initial datum $u_0$ satisfying
(cf.~\eqref{deficalVm})
\beeq{uzeroweak}
  u_0\in\calH_m:=\big\{v\in H:|v\OO|\le m\big\}.
\end{equation}
To show this, we proceed in a somehow reverse order,
by first deriving some estimates and then inferring a
precise statement. We notice that still a rigorous procedure
should rely on approximation and passage to the limit arguments
(i.e.~working on~\Pneen\ or some analogue of its 
and then letting $n\nearrow\infty$).
Nevertheless, for brevity (and since all works
similarly with the previous Section)
we prefer to consider here directly,
although formally, a limit solution $u$.
Thus, let us set
\beeq{defimu}
  \mu(s):=\int_0^s\frac{\dir}{b(r)},\qquad
   \muciapo(s):=\int_0^s\mu(r)\,\dir.
\end{equation}
Of course, by \eqref{hpb}, $\muciapo$ satisfies
\beeq{propmu}
  \frac1\mu s^2\le 2\muciapo(s)\le \frac1\alpha s^2
   \quad\perogni s\in\RR.
\end{equation}
Let us now perform an estimate of {\sl entropy}\/
type. Namely, let us test \eqref{eqne} by $\mu(u)$,
\eqref{eqne2} by $Bu$, and take the sum. Noting
that a couple of terms cancel, we infer
\beeq{entropy}
  2 \ddt \io\muciapo(u)
   +\|Bu\|^2
   +2\io W''(u)|\nabla u|^2
  \le\|f\|^2.
\end{equation}
Adding $\|u\|^2+2\lambda\|\nabla u\|^2$ to both hands sides
and using \eqref{newW} and the Poincar\'e-Wirtinger inequality
\eqref{PoWi}, we readily obtain
\beeq{entropy2}
  2 \ddt \io\muciapo(u)
   +\|u\|_{H^2(\Omega)}^2
   +2\eta\io \big(|u|^{p-2}|\nabla u|^2\big)
  \le c_3\big(1+\|\nabla u\|^2\big),
\end{equation}
where $c_3$ on the \rhs\ depends on $\lambda$,
$m$ (cf.~\eqref{uzeroweak}) and on the $H$-norm of $f$.

Once the initial datum $u_0\in\calH_m$ is given, 
owing to \eqref{propmu}, and noting that
\beeq{entropy3}
  c_3\big(1+\|\nabla u\|^2\big)
   \le\frac12\|u\|^2_{H^2(\Omega)}+c\|u\|^2+c,
\end{equation}
an application of Gronwall's Lemma in \eqref{entropy2} gives, 
for $T>0$, the a priori estimate (notice that we control
the full $V$-norm of $|u|^{(p-2)/2}u$ since we know
that $|u\OO|\le m$)
\beeq{stw11}
  \|u\|_{\LIH\cap\LDHD}
   +\||u|^{(p-2)/2}u\|_{\LDV} 
  \le \cc.
\end{equation}
Here and below the constants $\cc>0$ may depend on $T$
and $u_0$. Hence, by Sobolev's embeddings, it is also
\beeq{stw12}
  \|u\|_{L^{p}(0,T;L^{3p}(\Omega))}
  \le \cc.
\end{equation}
To show existence of an {\sl entropy}\/ solution, we 
have to see that \eqref{stw11} and \eqref{stw12} are enough to 
take the limit in equations \eqref{eqne}--\eqref{eqne2}
(recall we should work on some approximation, here).
By interpolation of Lebesgue spaces,
we actually have (here we just use that $p>2$)
\beeq{stw13}
  \|u\|_{L^{(3p-2)(p-1)/3(p-2)}(0,T;L^{2(p-1)}(\Omega))}
   \le \cc,
\end{equation}
whence, by \eqref{growW} 
and with the help of
a comparison of terms in \eqref{eqne2}, we have 
%
%
%
\beeq{stw14bis}
  \|W'(u)\|_{L^{(3p-2)/3(p-2)}(0,T;H)}
   +\|w\|_{L^{(3p-2)/3(p-2)}(0,T;H)}
   \le \cc.
\end{equation}
However, the equation \eqref{eqne} makes no sense in that
form as $u$ has only the above regularity. Nevertheless, 
taking $v\in H^3_{\bnn}(\Omega)$ 
(the (closed) subspace of $H^3(\Omega)$ consisting of functions
with $0$ normal derivative on $\de\Omega$), 
a.e.~in $(0,T)$ we can formally write
\beeq{perparti}
  \duav{B_uw,v}
   =\io b(u)\nabla w\cdot\nabla v
   =-\io w b'(u)\nabla u\cdot\nabla v
    -\io b(u) w \Delta v.
\end{equation}
Thus \eqref{eqne} has the weak correspondent
\beeq{eqnew}
  \duav{u_t,v}
    -\io w b'(u)\nabla u\cdot\nabla v
    -\io b(u) w \Delta v=0.
\end{equation}
Let us prove that \eqref{eqnew} does make sense
in the present regularity framework. Actually, 
noticing that, by \eqref{hpb},
\beeq{stw15}
  \|b(u)\|_{L^\infty(\Omega\times(0,T))}
   +\|b'(u)\|_{L^\infty(\Omega\times(0,T))}
   \le c,
\end{equation}
we readily obtain that
\beeq{stw16}
  \frac{|\duav{u_t,v}|}{\|v\|_{H^3(\Omega)}}
   \le c \|w\| \big(\|\nabla u\| + 1\big).
\end{equation}
Thus, observing that by \eqref{stw11} and interpolation
$\nabla u$ is bounded in $L^4(0,T;H)$, noting that for $p<6$ it
is $(3p-2)/3(p-2)>4/3$, taking the 
supremum w.r.t.~$v\in H^3_{\bnn}(\Omega)$
in \eqref{stw16}, and integrating the result over $(0,T)$,
we get 
\beeq{stw17}
  \|u_t\|_{L^{1}(0,T;(H^3_{\bnn})'(\Omega))}
  \le c_T,
\end{equation}
which shows that \eqref{eqnew} makes sense, provided that we interpret
the first term as a duality between $H^3_{\bnn}(\Omega)$ and its dual
(notice that, in fact, $H^3_0(\Omega)$ is contained
in $H^3_{\bnn}(\Omega)$, so that, in particular, \eqref{eqnew}
can be seen a relation in $H^{-3}(\Omega)$).

Let us now see that, more precisely,
\eqref{eqne} passes to the limit (i.e., that we have 
sufficient compactness to remove some kind of 
approximation). Actually, by \eqref{hpb}, \eqref{stw11},
\eqref{stw17} and the generalized Aubin Lemma
\cite[Cor.~4]{Si}, we have the convergences (holding
at least for suitable subsequences of the approximating
solutions, as usual) of
\beal{con11}
  & b(u),b'(u), \quext{strongly in }\/L^a(\Omega\times(0,\infty))~~
      \perogni a\in[1,\infty),\\
 \label{con12}
  & u, \quext{strongly in }\/L^{b}(0,T;V)~~\perogni b\in[1,4).
\end{align}
Moreover, using that $p$ is {\sl strictly}\/ lower than $6$ and
modifying a bit the argument leading to \eqref{stw14bis}, we 
can show that there exist exponents $a_*>4/3$ and $b_*>2$ such
that 
\beeq{con13}
  w \quext{converges weakly in }\/L^{a_*}(0,T;L^{b_*}(\Omega)).
\end{equation}
Thus, it is easy to see that \eqref{con11}--\eqref{con13}
allow us to pass to the limit in \eqref{eqnew}. 
Also the limit of \eqref{eqne2} is then 
easily taken since \eqref{con12} and
the weak convergence coming from the first of \eqref{stw14bis}
easily allow to identify the limit of $W'(u)$ by the usual
monotonicity argument. We have thus proven the 
\bete\label{teoesientro}
 Let\/ \eqref{hpb}, \eqref{W}, \eqref{growW}
 and\/ \eqref{newW} with $p\in(2,6)$ hold. Let $f$
 and $u_0$ satisfy\/ \eqref{regof} and\/ 
 \eqref{uzeroweak}, respectively. Then, there exists 
 at least one couple $(u,w)$ complying with the 
 regularity properties\/ \eqref{stw11}, \eqref{stw12},
 \eqref{stw14bis} and\/ \eqref{stw17}
 and such that\/ \eqref{eqnew} holds for all 
 $v\in H^3_{\bnn}(\Omega)$ and a.e.~in $(0,T)$.
 Moreover, it is 
 \beeq{eqne2w}
   (w-W'(u))=Bu+f \qquext{in }\/H,\quext{a.e.~in }\/(0,T),
 \end{equation}
 and the initial condition\/ \eqref{iniz} holds
 in $H$ (indeed, by \eqref{stw11} and \eqref{stw17},
 $u\in C_w([0,T];H)$). We call such a function $u$ an\/
 {\rm entropy solution} to\/ {\rm Problem (P$_0$)}.
\ente
Let us now study the long time behavior of
entropy solutions which, in a sense that might
be specified following the lines of the previous two 
Sections (cf.~Definition~\ref{defisollimi}),
have a {\sl limiting}\/ character. Still,
we prefer to work formally and do not enter the
details of the approximation-limit argument,
which should be very close to that
sketched in the previous
Sections. It is anyway worth noting
that for entropy solutions
(which are less regular than ``energy'' ones)
we do not expect, {\sl a fortiori}, any 
uniqueness property.
Our last result in this paper is the following
\bete\label{teoattraentro}
 Under the assumptions of\/ {\rm Theorem~\ref{teoesientro}},
 the set $\calSen$ of {\rm limiting} entropy solutions 
 to\/ {\rm (P$_0$)} constitutes a\/ {\rm limiting semiflow}
 on $\calH_m$. Furthermore, $\calSen$ admits 
 the global attractor $\calAen$ which is compact in 
 $\calH_m$ and bounded in $\calV_m$.
\ente
\noindent
{\bf Proof of Theorem~\ref{teoattraentro}}.~~%
The key point is to show that the entropy estimate
\eqref{entropy2} derived above also has a {\sl dissipative}\/ 
character. Let us then take $M>0$ 
(whose value will be chosen later) and set
\beeq{uM}
   u_M:=\max\big\{-M,\min\{u,M\}\big\}.
\end{equation}
Then, it is clear that
\beeq{con21}
  c_3\big(1+\|\nabla u\|^2\big)
   \le c_3+\frac{c_3}{M^{p-2}}
       \io\big(|u|^{p-2}|\nabla u|^2\big)
  +c_3\io|\nabla u_M|^2
\end{equation}
and, by interpolation, for all $\sigma>0$
the latter term can be controlled this way:
\begin{align}
 \nonumber
  c_3\io|\nabla u_M|^2
   & \le \sigma\|u_M\|_{H^{5/4}(\Omega)}^2
    +c(\sigma)\|u_M\|^2\\
 \label{con22}
   & \le c_4\sigma\|u\|_{H^{5/4}(\Omega)}^2
       +c(\sigma)\|u_M\|^2
     \le c_5\sigma\|u\|_{H^2(\Omega)}^2
       +c(\sigma,M,\Omega).
\end{align}
We used here the fact that the truncation operator
$u\mapsto u_M$ is continuous from $H^s(\Omega)$ to
itself for all $s<3/2$ (cf., e.g., \cite[Remark~0.1]{Sa}).

Thus, choosing $\sigma$ such that $c_5\sigma=1/2$ 
and $M$ so large that $c_3/M^{p-2}\le\eta$, 
\eqref{entropy2} gives 
\beeq{entropy2bis}
  2 \ddt \io\muciapo(u)
   +\frac12\|u\|_{H^2(\Omega)}^2
   +\eta\io \big(|u|^{p-2}|\nabla u|^2\big)
  \le c_6,
\end{equation}
where $c_6$ depends on $m$ but is independent of the 
choice of $u_0$. By \eqref{propmu} and
Gronwall's Lemma, \eqref{entropy2bis}
readily gives dissipativity in the space $\calH_0$ as
well as, for all $\tau>0$, the analogue
of \eqref{stunif1}, namely
\beeq{stunif1newbis}
  \sup_{t\in(\tau,\infty)}
   \int_t^{t+1}\|u(s)\|_{H^2(\Omega)}^2\,\dis
    \le c_7\|u_0\|^2e^{-\kappa'\tau}+C_0',
\end{equation}
for suitable $c_7,\kappa',C_0'$ independent of $u_0$. 
At this point, writing the energy equality in the 
form \eqref{energy} (with $\calE$ in place of $\calE_n$)
and using the uniform Gronwall Lemma, we immediately obtain 
existence of a uniformly absorbing set bounded in
$\calV$ and hence compact in $\calH$. This 
fact implies existence of the attractor
and concludes the proof.
\beos
 Coming back to the local compactness argument 
 in the previous Section, one can readily see
 that $\calAen$ is not only bounded, 
 but even {\sl compact}, in the space $\calV$.
\eddos




\vspace{10mm}

\noindent%
{\bf Author's address:}\\[1mm]
Giulio Schimperna\\
Dipartimento di Matematica, Universit\`a degli Studi di Pavia\\
Via Ferrata, 1,~~I-27100 Pavia,~~Italy\\
E-mail:~~{\tt giusch04@unipv.it}

\end{document}